\documentclass[titlepage,11pt]{article}
\usepackage{latexsym}
\usepackage{amsfonts}
\usepackage{amsmath}
\usepackage{tikz, fp}
\usepackage{textcomp}
\usetikzlibrary{graphs}
\usepackage{float}
\usetikzlibrary{decorations.markings}
\usetikzlibrary{decorations,decorations.pathmorphing}
\newcommand\blackslug{\hbox{\hskip 1pt \vrule width 4pt height 8pt depth 1.5pt
        \hskip 1pt}}
\newcommand\bbox{\hfill \quad \blackslug \bigbreak}

\def\LL{,\ldots,}
\def\cupcup{\cup\cdots\cup}
\newcommand{\vare}{\varepsilon}

\title{Pure pairs. IV. Trees in bipartite graphs}
\author{Alex Scott\\
Mathematical Institute, University of Oxford, Oxford OX2 6GG, UK
\\
\\
Paul Seymour\thanks{Supported by AFOSR grants A9550-19-1-0187 and FA9550-22-1-0234,  and NSF
grants DMS-1265563, DMS-1800053 and DMS-2154169.}\\
Princeton University, Princeton, NJ 08544
\\
\\
Sophie Spirkl\thanks{This material is based upon work supported by the National Science
Foundation under Award No. DMS-1802201. Current address: University of Waterloo, Waterloo, Ontario N2L3G1, Canada}\\
Princeton University, Princeton, NJ 08544.}

\date{October 14, 2019; revised \today}

\newtheorem{thm}{}[section]

\newcommand{\Proof}{\noindent{\bf Proof.}\ \ }

\begin{document}
\maketitle
\begin{abstract}
In this paper we investigate the bipartite analogue of the strong Erd\H os-Hajnal property.
We prove that for every forest $H$ and every $\tau$ with $0<\tau\le 1$, there exists
$\vare>0$, such that if $G$ has a bipartition $(A,B)$ and does not contain $H$ as an induced subgraph, and has at most 
$(1-\tau)|A|\cdot|B|$ edges, then 
there is a stable set $X$ of $G$ with $|X\cap A|\ge \vare |A|$ and $|X\cap B|\ge \vare |B|$.
No graphs $H$ except forests have this property. 

\end{abstract}

\section{Introduction}

All graphs in this paper are finite and have no loops
or parallel edges.
If $G,H$ are graphs, we say $G$ {\em contains} $H$ if some induced subgraph of $G$ is isomorphic to $H$, and $G$ is
{\em $H$-free} otherwise. 
We denote by $\alpha(G), \omega(G)$ the cardinalities of the largest stable sets and largest cliques in $G$
respectively. 
Two disjoint sets $A,B$ are {\em complete} if every vertex in $A$ is adjacent to every vertex in $B$, and
{\em anticomplete} if no vertex in $A$ has a neighbour in $B$; and we say $A$ {\em covers} $B$ if
every vertex in $B$ has a neighbour in $A$. A pair $(A,B)$ of subsets of $V(G)$ is {\em pure} if $A,B$ are 
complete or anticomplete to each other. We denote the complement graph of $H$ by $\overline{H}$. We denote the number of
vertices of $G$ by $|G|$.

The Erd\H{o}s-Hajnal conjecture~\cite{EH0, EH} asserts that:
\begin{thm}\label{EHconj}
{\bf Conjecture: }For every graph $H$, there exists $c>0$ such that every $H$-free graph $G$ satisfies
$$\alpha(G)\omega(G)\ge |G|^c.$$
\end{thm}
An {\em ideal} or {\em hereditary class} of graphs is a class $\mathcal{G}$ of 
graphs such that if $G\in \mathcal{G}$ and $H$ is isomorphic to an induced subgraph of $G$ then $H\in \mathcal{G}$.
We say that an ideal $\mathcal G$ of graphs has the {\em Erd\H os-Hajnal property} if
there exists $c>0$ such that every $G\in\mathcal G$ satisfies
$\alpha(G)\omega(G)\ge |G|^c$.  Thus the  Erd\H{o}s-Hajnal conjecture states that the class of $H$-free graphs has the Erd\H os-Hajnal property.

One way to prove that a class of graphs has the Erd\H os-Hajnal property is to prove something stronger.
We say that a class $\mathcal G$  has the {\em strong Erd\H os-Hajnal property} if
there is some $\vare>0$ such that every graph $G\in\mathcal G$ with at least two vertices 
contains disjoint sets $A$, $B$ that have size at least $\vare |G|$ and are either complete or anticomplete  (that is, $(A,B)$ is a pure pair).
It is easy to prove that, for an ideal $\mathcal G$, the strong Erd\H{o}s-Hajnal property implies the Erd\H{o}s-Hajnal property 
(see~\cite{APPRS,fp}).  Unfortunately, it is also easy to prove (by considering a sparse random graph with girth larger than $|H|$) that if the class of $H$-free graphs has the 
strong Erd\H os-Hajnal property then $H$ must be a forest; and (by considering the complement of a sparse random graph) that $\overline H$ must also be a forest.  Thus $H$ has at most four vertices, and the Erd\H os-Hajnal conjecture was already known for these graphs.

But what if we exclude {\em both} a forest and the complement of a forest?
Then there is some good news.   In an earlier paper~\cite{trees}, with Maria Chudnovsky, we proved that this implies the strong Erd\H os-Hajnal property:
\begin{thm}\label{oldforestsymm}
For every forest $H$, there exists $\vare>0$ such that for every graph $G$ that is both $H$-free and $\overline{H}$-free 
with $|G|\ge 2$, there is a pure pair $(Z_1,Z_2)$ of subsets of $V(G)$ with $|Z_1|,|Z_2|\ge \vare |G|$.
\end{thm}
We also proved the stronger result that, if $G$ is not too dense, then it is enough just to exclude a forest.
\begin{thm}\label{oldmainthm}
For every forest $H$ there exists $\vare>0$ such that for every $H$-free graph $G$ with $|G|\ge 2$, either
\begin{itemize}
\item some vertex has degree at least $\vare |G|$; or
\item there exist disjoint $Z_1,Z_2\subseteq V(G)$ with $|Z_1|,|Z_2|\ge \vare |G|$, anticomplete.
\end{itemize}
\end{thm}
Neither of these theorems hold for any graph $H$ such that neither of $H, \overline{H}$ is a forest
(this follows easily from the random construction 
by Erd\H{o}s of graphs with large girth and 
large chromatic number~\cite{girth}).

In this paper we look at the analogous question for bipartite graphs. 
It will help to 
set up some terminology. 
A {\em bigraph} $G$ is a graph with a designated bipartition $(V_1(G), V_2(G))$; thus, $V_1(G), V_2(G)$ are disjoint stable sets of $G$ with union $V(G)$.  The bigraph obtained from the same graph by
exchanging $V_1(G)$ and $V_2(G)$ is called the {\em transpose} of $G$.
The {\em bicomplement} of a bigraph $G$ is the bigraph $G'$ with the same vertex set, and the same bipartition,
in which for all $v\in V_1(G)$ and $w\in V_2(G)$, $v$ and $w$ are adjacent in $G$ if and only if they are not adjacent in $G'$.
If $G$ is a bigraph, and $X\subseteq V(G)$, $G[X]$ denotes the bigraph induced on $X$ in the natural sense (that is, the bigraph 
$H$ with $V(H)=X$, where for all $u,v\in X$, $u,v$ are adjacent in $H$ if and only if they are adjacent in $G$, and 
$V_i(H)=V_i(G)\cap X$ for $i = 1,2$). We call $G[X]$ an {\em induced sub-bigraph} of $G$. An {\em isomorphism} between 
bigraphs $G,H$ is an isomorphism between the corresponding graphs
that maps $V_i(G)$ to $V_i(H)$ for $i = 1,2$.
If $G,H$ are bigraphs, we say that $G$ {\em contains} $H$ if there is an isomorphism from $H$ to an induced sub-bigraph of $G$.
We say $G$ is {\em $H$-free} if $G$ does not contain $H$. We stress that a bigraph $G$ may contain a bigraph $H$ and not contain 
the transpose of $H$.

What about the analogue of \ref{oldforestsymm} for bigraphs? Let $H$ be a forest bigraph, and suppose that $G$ is a bigraph that 
contains neither $H$ not its bicomplement: 
must there then be a pure pair of subsets of $V_1(G)$ and $V_2(G)$ respectively, both of linear size? 
This seems neither to imply, nor to be implied by, the results of~\cite{trees} just mentioned,
but there has been some previous work on it. It was 
conjectured in this form by Alecu, Atminas, Lozin and Zamaraev~\cite{Alecu},
amd Axenovich, Tompkins and Weber~\cite{axenovich}, 
and independently Kor\'andi, Pach and Tomon~\cite{pach} gave a result for bigraphs $H$ in which $|V_1(H)|\le 2$:
\begin{thm}\label{twovert}
Let $H$ be a forest bigraph with $|V_1(H)|\le 2$, such that the bicomplement of $H$ is also a forest.
Then there exists $\vare>0$ such that, if $G$ is a bigraph that is $H$-free then
there is a pure pair $(Z_1,Z_2)$ in $G$ with $Z_i\subseteq V_i(G)$ and
$|Z_i|\ge \vare |V_i(G)|$ for $i = 1,2$.
\end{thm}
In this paper we will prove the conjecture of Alecu, Atminas, Lozin and Zamaraev~\cite{Alecu},
that the bipartite analogue of \ref{oldforestsymm} holds in full. Our main result is the following:
\begin{thm}\label{forestsymm}
For every forest bigraph $H$ with bicomplement $J$, there exists $\vare>0$ such that, if $G$ is a bigraph that is 
both $H$-free and $J$-free, then there is a pure pair $(Z_1,Z_2)$ in $G$ with $Z_i\subseteq V_i(G)$ and
$|Z_i|\ge \vare |V_i(G)|$ for $i = 1,2$.
\end{thm}

Note that this implies that it is sufficient to exclude any forest $H$ and any forest bicomplement $J$ (as we can always consider the forest obtained from the disjoint union of $H$ and $\overline J$).  Furthermore, if we exclude any finite set of bigraphs, then the random construction once again shows that one of the graphs must be a forest, and one must be the bicomplement of a forest.   Thus \ref{forestsymm} leads to a characterizion of the finite sets of excluded subgraphs that give the (bipartite) strong Erd\H os-Hajnal property.

As with \ref{oldmainthm}, we will also prove a `one-sided' version of the result for sparse graphs, where we only exclude a forest.
Let us say a bigraph $G$ is {\em $\vare$-coherent}, where $\vare>0$, if:
\begin{itemize}
\item every vertex in $V_1(G)$ has degree less than $\vare |V_2(G)|$;
\item every vertex in $V_2(G)$ has degree less than $\vare |V_1(G)|$; 
\item there do not exist anticomplete subsets $Z_1\subseteq V_1(G)$ and $Z_2\subseteq V_2(G)$, such that $|Z_i|\ge \vare |V_i(G)|$ for $i = 1,2$.
\end{itemize}
Thus, if $0<\vare\le \vare'$ and $G$ is $\vare$-coherent then it is also $\vare'$-coherent.
Our main theorem says:

\begin{thm}\label{mainthm}
For every forest bigraph $H$, there exists $\vare>0$ such that every $\vare$-coherent bigraph contains $H$.
\end{thm}

This can be further strengthened, to prove the statement in the abstract; in the next section we will show that the ``sparse'' hypothesis 
can be replaced
with a ``not very dense'' hypothesis,
as follows:
\begin{thm}\label{betterthm}
For every forest bigraph $H$, and every $\tau$ with $0<\tau\le 1$, there exists $\vare>0$ such that if $G$ is an $H$-free bigraph with at
most $(1-\tau)|V_1(G)|\cdot|V_2(G)|$ edges, then there are anticomplete sets $Z_i\subseteq V_i(G)$ with $|Z_i|\ge \vare |V_i(G)|$
for $i = 1,2$.
\end{thm}

Finally, let us note that there are very interesting questions of this type for {\em ordered} bigraphs, or equivalently 0-1 matrices.
Kor\'andi, Pach and Tomon~\cite{pach} also made a much stronger conjecture, that the statement of \ref{betterthm} still holds 
if we work with ordered bigraphs (bigraphs $G$ with
fixed linear orders on $V_1(G)$ and on $V_2(G)$) and ask for containment that respects the orders. 
This remains open, and is discussed further in~\cite{orderedbiptrees}.

\section{Reducing to the sparse case}

In this section we will show that \ref{mainthm} implies \ref{betterthm} and \ref{forestsymm}.
We need the following. For general graphs an analogue was proved by R\"odl~\cite{rodl}, using the regularity lemma, but
for bigraphs its proof is much easier, and does not use the regularity lemma. The result is essentially due
to Erd\H{o}s, Hajnal and Pach~\cite{EHP}, but we give a proof since it is short.
\begin{thm}\label{getsparse}
Let $H$ be a bigraph, and let $\vare>0$. Then there exists $\delta>0$ with the following property. Let $G$ be an $H$-free bigraph 
with $V_1(G), V_2(G)\ne \emptyset$;
then there exists $A_i\subseteq V_i(G)$ with $|A_i|\ge \delta|V_i(G)|$ for $i = 1,2$, such that either
\begin{itemize}
\item every vertex in $A_1$ has fewer than $\vare|A_2|$ neighbours in $A_2$, and every vertex in $A_2$ has fewer than  
$\vare|A_1|$
neighbours in $A_1$; or
\item every vertex in $A_1$ has more than $(1-\vare)|A_2|$ neighbours in $A_2$, and every vertex in $A_2$ has more 
than $(1-\vare)|A_1|$
neighbours in $A_1$.
\end{itemize}
\end{thm}
\Proof We may assume that $\vare<1$ and $V(H)\ne \emptyset$, because otherwise the result is trivial. 
Let $V_1(H)=\{a_1\LL a_k\}$ and $V_2(H)=\{b_1\LL b_{\ell}\}$. If $k=0$ or $\ell=0$ the theorem holds with 
$\delta= \min(1/2,|H|^{-1})$, so we may assume that $k,\ell>0$. Define $\delta=\min(1/2,|H|^{-1},(\vare/2)^k/\ell)$.
We claim that $\delta$ satisfies the theorem. Let $G$ be an $H$-free bigraph.
If $|V_1(G)|< \delta^{-1}$, the result holds with 
$|A_1|=1, A_1=\{v\}$ say, and and $A_2$ the larger of $N(v)$ and $V_2(G)\setminus N(v)$, since $\delta\le 1/2$; so we may assume that 
$|V_1(G)|\ge  \delta^{-1}$, and similarly $|V_2(G)|\ge  \delta^{-1}$.

Let $0\le i\le k$, let $u_1\LL u_i\in V_1(G)$
be distinct, and let $1\le j\le \ell$. We say that $v\in V_2(G)$ is {\em $j$-appropriate} for $(u_1\LL u_i)$ if for all $h$ with $1\le h\le i$:
\begin{itemize}
\item if $a_h, b_j$ are adjacent in $H$ then $u_h,v$ are adjacent in $G$; and
\item if $a_h,b_j$ are nonadjacent in $H$ then $u_h,v$ are nonadjacent in $G$.
\end{itemize}
For all $i$ with $0\le i\le k$ let $n_i=(\vare/2)^i|V_2(G)|$. Thus 
$n_k=(\vare/2)^k|V_2(G)|\ge (\vare/2)^k\delta^{-1}\ge \ell$.
For $0\le i\le k$, we say a sequence $(u_1\LL u_i)$ of vertices in $V_1(G)$ is {\em $i$-good} if $u_1\LL u_i$ are all distinct, and
for $1\le j\le \ell$ there are at least $n_i$ vertices in $V_2(G)$ that are $j$-appropriate for $(u_1\LL u_i)$.

The null sequence is $0$-good, since  $n_0=|V_2(G)|$; and so we may
choose $i$ with $0\le i\le k$ maximum such that there is an $i$-good sequence  $(u_1\LL u_i)$.
Suppose that $i=k$.
For $1\le j\le \ell$ choose $v_j\in V_2(G)$, $j$-appropriate for $(u_1\LL u_k)$,
such that $v_1\LL v_{\ell}$ are all different (the last is possible since $n_k\ge \ell$); then the subgraph 
induced on $\{u_1\LL u_k, v_1\LL v_{\ell}\}$ is isomorphic to $H$,
a contradiction. So $i<k$.

From the maximality of $i$, for each $u\in V_1(G)\setminus \{u_1\LL u_i\}$, there exists $j\in \{1\LL \ell\}$ such that fewer than $n_{i+1}$
vertices in $V_2(G)$ are $j$-appropriate for $(u_1\LL u_i, u)$. Hence there exist $A_1\subseteq V_1(G)\setminus \{u_1\LL u_i\}$
with $|A_1|\ge (|V_1(G)|-k)/\ell$, and $j\in \{1\LL \ell\}$, such that for each $u\in A_1$, 
fewer than $n_{i+1}$
vertices in $V_2(G)$ are $j$-appropriate for $(u_1\LL u_i, u)$.
Let $A_2'$ be the set of vertices in $V_2(G)$ that are $j$-appropriate for $(u_1\LL u_i)$; thus $|A_2'|\ge n_i$, from the choice of $i$.
Consequently 
\begin{itemize}
\item if $a_{i+1}, b_j$ are adjacent in $H$, then each $u\in A_1$ has fewer than $n_{i+1}$ neighbours in $A_2'$;
\item if $a_{i+1}, b_j$ are nonadjacent in $H$, then each $u\in A_1$ has more than $|A_2'|-n_{i+1}$ neighbours in $A_2'$.
\end{itemize}
By taking bicomplements if necessary, we may assume the former. There are fewer than $n_{i+1}|A_1|$ edges between $A_1,A_2'$, and so
at least $|A_2'|/2$ vertices in $A_2'$ have fewer than $2n_{i+1}|A_1|/|A_2'|$ neighbours in $A_1$. Let $A_2$ be the set of vertices in $A_2'$ 
with fewer than $\vare |A_1|$ neighbours in $A_1$. Since $\vare \ge 2n_{i+1}/|A_2'|$ (because $|A_2'|\ge n_i=(2/\vare)n_{i+1}$), it follows
that $|A_2|\ge |A_2'|/2\ge n_i/2$. Thus:
\begin{itemize}
\item $|A_1|\ge \delta|V_1(G)|$, since $|A_1|\ge (|V_1(G)|-k)/\ell\ge \vare|V_1(G)|$ (because $|V_1(G)|\ge \vare^{-1}\ge k+\ell$);
\item $|A_2|\ge \delta|V_2(G)|$, since 
$|A_2| \ge n_i/2\ge n_k= (\vare/2)^k|V_2(G)|\ge \delta|V_2(G)|$;
\item every vertex in $A_1$ has fewer than $n_{i+1}= \vare n_i/2\le  \vare|A_2|$ neighbours in $A_2$; and
\item every vertex in $A_2$ has fewer than $\vare|A_1|$ neighbours in $A_1$.
\end{itemize}
This proves \ref{getsparse}.~\bbox

\noindent{\bf Proof of \ref{forestsymm}, assuming \ref{mainthm}.\ \ }
Let $H$ be a forest bigraph with bicomplement $J$. By \ref{mainthm} there exists $\vare>0$ such that
every $\vare$-coherent bigraph contains $H$. By \ref{getsparse} there exists $\delta'>0$ such that taking $\delta = \delta'$ satisfies
\ref{getsparse}. Let $\delta=\vare\delta'$; we claim that $\delta$ satisfies \ref{forestsymm}. 

Let $G$ be a bigraph that is both $H$-free and $J$-free. We must show that there exist $B_1\subseteq V_1(G)$ and 
$B_2\subseteq V_2(G)$ such that $|B_i|\ge \delta|V_i(G)|$ for $i = 1,2$, and $B_1, B_2$ are complete or anticomplete. 
We may assume that $V_i(G)\ne \emptyset$ for $i = 1,2$.
By the choice of $\delta'$, there exist $A_i\subseteq V_1(G)$ for $i = 1,2$, such that $|A_i|\ge \delta'|V_i(G)|$ for $i = 1,2$, and either
\begin{itemize}
\item every vertex in $A_1$ has fewer than $\vare|A_2|$ neighbours in $A_2$, and every vertex in $A_2$ has fewer than $\vare|A_1|$
neighbours in $A_1$; or
\item every vertex in $A_1$ has more than $(1-\vare)|A_2|$ neighbours in $A_2$, and every vertex in $A_2$ has more than $(1-\vare)|A_1|$
neighbours in $A_1$.
\end{itemize}
Suppose that the first holds.
By applying \ref{mainthm} to the subgraph $G'$ of $G$ induced on $A_1\cup A_2$, we deduce that $G'$ is not
$\vare$-coherent; and so 
there exist anticomplete subsets $B_1\subseteq A_1$ and $B_2\subseteq A_2$, such that $|B_i|\ge \vare |A_i|$ for $i = 1,2$.
But $\vare|A_i|\ge \vare\delta'|V_i(G)| = \delta|V_i(G)|$, as required.

If the second holds, then we apply \ref{mainthm} to the bicomplement $G'$ of the subgraph of $G$ induced on $A_1\cup A_2$,
and deduce that $G'$ is not $\vare$-coherent (since it is $H$-free, because $G$ is $J$-free). So 
there exist subsets $B_1\subseteq A_1$ and $B_2\subseteq A_2$, such that $|B_i|\ge \vare |A_i|$ for $i = 1,2$, and
$B_1, B_2$ are anticomplete in $G'$ and hence complete in $G$. Since $|B_i|\ge \delta|V_i(G)|$ for $i = 1,2$, the result follows.
This proves \ref{forestsymm}.~\bbox

We can in fact weaken the hypothesis of \ref{mainthm} that $G$ is sparse. 
In \ref{betterthm} we replace this by
the hypothesis that $G$ is not very dense, that the bicomplement of $G$
has at least $\tau n^2$ edges.
The ``not very dense'' hypothesis is as good as
the ``sparse'' hypothesis because of the next result.

\begin{thm}\label{countnonedges}
For all $c,\vare,\tau>0$ with $\vare<\tau\le 8/9$, there exists $\delta>0$ with the following property. Let $G$
be a bigraph with at most $(1-\tau)|V_1(G)|\cdot|V_2(G)|$  edges and with $V_1(G), V_2(G)\ne \emptyset$. Then
there exist $Z_i\subseteq V_i(G)$ with  $|Z_i|\ge \delta|V_i(G)|$ for $i = 1,2$, such that
there are fewer than
$(1-\vare)|Y_1|\cdot|Y_2|$ edges between $Y_1,Y_2$ for all subsets $Y_i\subseteq Z_i$ with
$|Y_i|\ge  c|Z_i|$ for $i = 1,2$.
\end{thm}
\Proof
By reducing $c$, we may assume that $c\le 1/3$. 
Let
$$\lambda= 1-\frac{(\tau-\vare)c^2}{(1-c^2)(1-\tau)}.$$
It follows that $0\le \lambda<1$ (since $c^2/(1-c^2)\le 1/8$ and $(\tau-\vare)/(1-\tau)< 8$).
Choose an integer $n\ge 0$ such that $\lambda^n(1-\tau)\le (1-\vare)c/2$. Let $\delta=\min(c^n, \tau)$. We will show that $\delta$
satisfies the theorem.

Let $G$
be a bigraph with at most $(1-\tau)|V_1(G)|\cdot|V_2(G)|$  edges.
Choose an integer $t\ge 0$ with $t\le n$, maximum such that there are subsets $Z_i\subseteq V_i(G)$ with $|Z_i|\ge c^t|V_i(G)|$
for $i = 1,2$, where the number of edges
between $Z_1,Z_2$ is at most $\lambda^t(1-\tau)|Z_1|\cdot|Z_2|$.
(This is possible since we may take $t=0$ and $Z_i=V_i(G)$ for $i = 1,2$.)
\\
\\
(1) {\em If $t=n$ then the theorem holds.}
\\
\\
Suppose that $t=n$. Thus there are at most $\lambda^n(1-\tau)|Z_1|\cdot|Z_2|\le (1-\vare)(c/2)|Z_1|\cdot|Z_2|$ edges
between $Z_1,Z_2$.
At least half of the vertices in $Z_1$ have at most $(1-\vare)c |Z_2|$ neighbours in $Z_2$; choose $Z_1'\subseteq Z_1$
with $|Z_1'|\ge |Z_1|/2\ge \delta|V_1(G)|$ such that every vertex in $Z_1'$ has at most $(1-\vare)c|Z_2|$ neighbours in $Z_2$.
Now let $Y_1\subseteq Z_1'$, and let $Y_2\subseteq Z_2$ with $|Y_2|\ge c|Z_2|$. Each vertex in
$Y_1$ has at most $(1-\vare)c|Z_2|\le (1-\vare) |Y_2|$ neighbours in $Y_2$, and so the number of edges between $Y_1,Y_2$ is at most
$(1-\vare)|Y_1|\cdot |Y_2|$. Since $|Z_1'|\ge \delta|V_1(G)|$ and $|Z_2|\ge \delta|V_2(G)|$, the pair $Z_1', Z_2$ satisfies the theorem.
This proves (1).
\\
\\
(2) {\em If $|V_1(G)|\le 1/\delta$ or $|V_2(G)|\le 1/\delta$ then the theorem holds.}
\\
\\
Suppose that $|V_1(G)|\le 1/\delta$, say.
Since $G$ has at most $(1-\tau)|V_1(G)|\cdot|V_2(G)|$ edges, and $V_1(G)\ne \emptyset$, some vertex $v_1\in V_1(G)$
has at most $(1-\tau)|V_2(G)|$ neighbours in $V_2(G)$; and so there is a set $Z_2\subseteq V_2(G)$ with
$|Z_2|\ge \tau|V_2(G)|\ge \delta|V_2(G)|$
that is anticomplete to $Z_1=\{v_1\}$. Since $|Z_1| =1\ge \delta |V_1(G)|$, the theorem holds. This proves (2).
\\
\\
(3) {\em If $t<n$ then the theorem holds.}
\\
\\
Suppose that $t<n$. Since 
$$|Z_i|\ge c^t|V_i(G)|\ge c^{n-1}|V_i(G)|\ge (\delta/c)|V_i(G)|\ge  \delta|V_i(G)|$$ 
for $i=1,2$, it suffices to show that there
do not exist subsets $Y_i\subseteq Z_i$ with $|Y_i|\ge c|Z_i|$ for $i = 1,2$, such that
the number of edges between $Y_1,Y_2$ is at least $(1-\vare)|Y_1|\cdot|Y_2|$.
Suppose that such subsets exist. By averaging, we may assume that $|Y_i|=\lceil c|Z_i|\rceil$ for $i = 1,2$;
by (2) we may assume that $\delta|V_i(G)|>1$ for $i = 1,2$, and therefore
$|Z_i|\ge (\delta/c) |V_i(G)|>1/c\ge 3$; and since $c\le 1/3$,  it follows that 
$$\lceil c|Z_i|\rceil\le c|Z_i|+1\le (1-c)|Z_i|.$$
Let $X_i=Z_i\setminus Y_i$, and let
$y_i=|Y_i|$ and $x_i=|X_i|$ for $i = 1,2$. Thus $x_i, y_i\ge c|Z_i|\ge c^{t+1}|V_i(G)|$, for $i = 1,2$. From the choice of $t$,
it follows that
\begin{itemize}
\item there are more than $\lambda^{t+1}(1-\tau)y_1x_2$ edges between $Y_1$ and $X_2$;
\item
there are more than $\lambda^{t+1}(1-\tau)x_1y_2$ edges between $X_1$ and $Y_2$; and
\item there are more than $\lambda^{t+1}(1-\tau)x_1x_2$ edges between $X_1$ and $X_2$.
\end{itemize}
Adding, we deduce that there are more than
$$\lambda^{t+1}(1-\tau)(y_1x_2+x_1y_2+x_1x_2) + (1-\vare)y_1y_2$$
edges between $Z_1,Z_2$, and so
$$\lambda^{t+1}(1-\tau)(y_1x_2+x_1y_2+x_1x_2) + (1-\vare)y_1y_2< \lambda^t(1-\tau)|Z_1|\cdot|Z_2|.$$
Since $|Z_i|=x_i+y_i$ for $i = 1,2$, it follows that
$$(1-\vare-\lambda^{t+1}(1-\tau))y_1y_2< \lambda^t(1-\lambda)(1-\tau)(x_1+y_1)(x_2+y_2),$$
and since $y_1y_2\ge c^2(x_1+y_1)(x_2+y_2)$, we deduce that
$$(1-\vare-\lambda^{t+1}(1-\tau))c^2< \lambda^t(1-\lambda)(1-\tau).$$
Since $\lambda^{t+1}\le \lambda$, and $\lambda^t\le 1$, it follows that
$(1-\vare-\lambda(1-\tau))c^2< (1-\lambda)(1-\tau),$
and so
$$\lambda(1-\tau)(1-c^2)< (1-\tau)-(1-\vare)c^2,$$
contradicting the definition of $\lambda$. This proves (3).

\bigskip

From (1) and (3), this proves \ref{countnonedges}.~\bbox

Let us deduce \ref{betterthm}, which we restate:
\begin{thm}\label{betterthm2}
For every forest bigraph $H$, and every $\tau>0$, there exists $\vare>0$ such that if $G$ is an $H$-free bigraph with at
most $(1-\tau)|V_1(G)|\cdot|V_2(G)|$ edges, then there are anticomplete sets $Z_i\subseteq V_i(G)$ with $|Z_i|\ge \vare |V_i(G)|$
for $i = 1,2$.
\end{thm}
\noindent{\bf Proof, assuming \ref{mainthm}.\ \ }
Let $H$ be a forest bigraph, and let $\tau>0$. By reducing $\tau$, we may assume that $\tau\le 8/9$.
By \ref{mainthm} there exists $\eta>0$ such that 
every $\eta$-coherent bigraph contains $H$. By reducing $\eta$, we may assume that $\eta<\tau$.
Choose $c>0$ such that setting $\delta=c$ satisfies \ref{getsparse} 
with $\vare$ replaced by $\eta$. Let $\delta$ satisfy \ref{countnonedges} with $\vare$ replaced by $\eta$.
Let $\vare=c\delta\eta$.

Now let $G$ be an $H$-free bigraph with at
most $(1-\tau)|V_1(G)|\cdot|V_2(G)|$ edges. We may assume that $V_1(G), V_2(G)\ne \emptyset$.
By \ref{countnonedges}, 
there exist $Z_i\subseteq V_i(G)$ with  $|Z_i|\ge \delta|V_i(G)|$ for $i = 1,2$, such that
there are fewer than
$(1-\eta)|Y_1|\cdot|Y_2|$ edges between $Y_1,Y_2$ for all subsets $Y_i\subseteq Z_i$ with
$|Y_i|\ge  c|Z_i|$ for $i = 1,2$. By \ref{getsparse}, applied to the sub-bigraph induced on $Z_1\cup Z_2$,
there exists $Y_i\subseteq Z_i$ with $|Y_i|\ge c|Z_i|$ for $i = 1,2$, such that either
\begin{itemize}
\item every vertex in $Y_1$ has fewer than $\eta|Y_2|$ neighbours in $Y_2$, and every vertex in $Y_2$ has fewer than
$\eta|Y_1|$
neighbours in $Y_1$; or
\item every vertex in $Y_1$ has more than $(1-\eta)|Y_2|$ neighbours in $Y_2$, and every vertex in $Y_2$ has more
than $(1-\eta)|Y_1|$
neighbours in $Y_1$.
\end{itemize}
The second is impossible, since for all such $Y_1,Y_2$ there are fewer than $(1-\eta)|Y_1|\cdot|Y_2|$ edges between $Y_1,Y_2$.
Thus the first bullet holds. Since the sub-bigraph induced on $Y_1\cup Y_2$ is $H$-free, it is not $\eta$-coherent; and
so there exist anticomplete sets $X_i\subseteq Y_i$ with 
$$|X_i|\ge \eta|Y_i|\ge \eta (c\delta|V_i(G)|)\ge \vare|V_i(G)|$$ 
for $i = 1,2$. 
This proves \ref{betterthm2}.~\bbox

\section{Parades and concavity}

In this section we carry out the main step in the proof of \ref{mainthm}.  The general approach is similar to the main proof in~\cite{trees}, and we apologize for repeating some material and ideas from there. But adapting the proof of~\cite{trees} to work for bipartite graphs was nontrivial, and there seems no way
to present the tricky new parts of the argument without also including the straightforward parts.

We will need a number of definitions.
A {\em parade} in a bigraph $G$ means a sequence 
$$(A_1\LL A_K; B_1\LL B_L)$$ 
of pairwise disjoint nonempty subsets of $V(G)$, 
such that 
\begin{itemize}
\item $A_1\LL A_K\subseteq V_1(G)$, and $B_1\LL B_L\subseteq V_2(G)$;
\item $A_1\LL A_K$ all have the same cardinality, and $B_1\LL B_L$ all have the same cardinality (possibly different).
\end{itemize}
Its {\em length} is the pair $(K,L)$,
and its {\em width} is the pair $(|A_1|, |B_1|)$.  (For convenience in handling widths, let us $(W_1',W_2')\le (W_1,W_2)$ if $W_i'\le W_i$ for $i = 1,2$, and define
$\lambda(W_1,W_2)=(\lambda W_1,\lambda W_2)$ for $\lambda\ge 0$.)
We call the sets $A_i, B_i$ {\em blocks} of the parade.
We are interested in parades of some fixed length, and width at least linear in $(|V_1(G)|,|V_2(G)|)$. (We remark that in later papers of this series we use ``parade'' to mean the same thing with the second bullet above removed; but here
it is convenient to include the second bullet in the definition.)

Here are two useful ways to make smaller parades from larger.
Let $\mathcal{P}=(A_1\LL A_K; B_1\LL B_L)$ be a parade. First, let $1\le r_1<r_2\cdots <r_k\le K$, and $1\le s_1<\cdots <s_{\ell}\le L$; 
then
$\mathcal{P}'=(A_{r_1}\LL A_{r_k}; B_{s_1}\LL B_{s_{\ell}})$ is a parade, of smaller length
but of the same width, and we call it a {\em sub-parade} of $\mathcal{P}$.
Let $I=\{r_1\LL r_k\}$ and $J=\{s_1\LL s_{\ell}\}$; then $\mathcal{P}[I;J]=(A_i\;(i\in I);B_j\;(j\in J))$ denotes the same subparade 
$\mathcal{P}'$.
Second, for $1\le i\le K$ let $A_i'\subseteq A_i$, all of the same cardinality, and for $1\le j\le L$ let $B_{j}'\subseteq B_{j}$, all of
the same cardinality; then
the sequence $(A_1'\LL A_K';B'_1\LL B'_L)$ is a parade, of the same length but of smaller width, and we call it
a {\em contraction} of $\mathcal{P}$. A contraction of a sub-parade (or equivalently, a sub-parade of a contraction)
we call a {\em minor} of $\mathcal{P}$. (Thus a minor of a minor is a minor.)

Let $\mathcal{P}=(A_1\LL A_K; B_1\LL B_L)$ be a parade in a bigraph $G$. We say an induced sub-bigraph $H$ of $G$ is {\em $\mathcal{P}$-rainbow} or {\em rainbow relative to $\mathcal{P}$}
if each vertex of $H$ belongs to some block of $\mathcal{P}$, and no two vertices belong to the same block.
A {\em copy} of a bigraph $T$ in a bigraph $G$ is a bigraph isomorphic to $T$ that is
contained in $G$.

An {\em ordered bigraph} is a bigraph $T$ with linear orders imposed on $V_1(T)$ and on $V_2(T)$.
Let $\mathcal{P}=(A_1\LL A_K; B_1\LL B_L)$ be a parade in $G$. If $H$ is an $\mathcal{P}$-rainbow
induced sub-bigraph of $G$, then there is an associated linear order $<$ on $V_1(H)$ defined by $u<v$ if $u\in A_i$ and $v\in A_j$
for some $i,j$ with $i<j$; and similarly for $V_2(H)$. This gives an ordered bigraph that we call the {\em $\mathcal{P}$-ordering}
of $H$; and if the $\mathcal{P}$-ordering of $H$ is isomorphic to some ordered bigraph $T$ we say that $H$ is a {\em copy} of $T$.

A {\em rooted bigraph} $H$ is a pair $(H^-,r(H))$, where $H^-$ is a bigraph and $r(H)\in V(H^-)$; we call $r(H)$ the {\em root}.
Thus, the root might belong to $V_1(H)$ or to $V_2(H)$.
If $H_1,H_2$ are rooted bigraphs, by an {\em isomorphism} between them we mean an isomorphism
between $H_1^-$ and $H_2^-$ that takes root to root.

An induced rooted sub-bigraph $H$ of $G$ is {\em $\mathcal{P}$-left-rainbow}
if 
\begin{itemize}
\item it is $\mathcal{P}$-rainbow;
and 
\item if the root of $H$ belongs to $A_h$, then $h\le i$ for all $i\in \{1\LL K\}$ with $V(H)\cap A_i\ne \emptyset$; and
if the root of $H$ belongs to $B_h$, then $h\le j$ for all $j\in \{1\LL L\}$ with $V(H)\cap B_j\ne \emptyset$.
\end{itemize}
We define {\em $\mathcal{P}$-right-rainbow} similarly, requiring $h\ge i$ and $h\ge j$ instead.

Let $\mathcal{P}=(A_1\LL A_K; B_1\LL B_L)$ be a parade in a bigraph $G$. If $T$ is an induced subgraph that is $\mathcal{P}$-rainbow,
its {\em support} is the pair $(I,J)$ of subsets of $\{1\LL K\}$ where
\begin{itemize}
\item $I$ is the set of all $i\in \{1\LL K\}$ such that $V(T)\cap A_i\ne \emptyset$, and 
\item $J$ is the set of all $j\in \{1\LL L\}$ such that $V(T)\cap B_j\ne \emptyset$.
\end{itemize}
If $S$ is an ordered bigraph, and $\mathcal{P}$ is as above, 
we define the {\em trace} of $S$ (relative to $\mathcal{P}$) to be the set of supports of all 
$\mathcal{P}$-rainbow copies of $S$ in $G$.

We say $\mathcal{P}$ is {\em $\tau$-support-uniform} if
for every ordered tree bigraph $T$ with at most $\tau$ vertices, either the trace of $T$ (relative to $\mathcal{P}$)
is empty, or it consists of all pairs $(I,J)$ with $I\subseteq \{1\LL K\}$ and $J\subseteq \{1\LL L\}$ of cardinalities $|V_1(T)|, |V_2(T)|$ respectively.

Let $0< \kappa\le 1$, 
and let $\mathcal{P}=(A_1\LL A_K; B_1\LL B_L)$ be a parade in a bigraph $G$. We say that $\mathcal{P}$ is 
{\em $(\kappa,\tau)$-support-invariant} if it has the following property:
for every contraction $\mathcal{P}'=(A_1'\LL A_K';B_1'\LL B_L')$ of $\mathcal{P}$ such that $|A_i'|\ge \kappa |A_i|$ for $1\le i\le K$
and $|B_j'|\ge \kappa|B_j|$
for $1\le j\le L$, and 
for every ordered tree bigraph $T$ with at most $\tau$ vertices,
the trace of $T$ relative to $\mathcal{P}$ equals the trace of $T$ relative to $\mathcal{P'}$.

If $G$ is a bigraph and $X\subseteq V_1(G)$ and $Y\subseteq V_2(G)$, or vice versa, and $0\le \lambda\le  1$, 
we say that $X$ {\em $\lambda$-covers}
$Y$ if there are at least $\lambda|Y|$ vertices in $Y$ with a neighbour in $X$, and $X$ {\em $\lambda$-misses}
$Y$ if there are at least $\lambda|Y|$ vertices in $Y$ with no neighbour in $X$.

A parade $\mathcal{P}=(A_1\LL A_K;B_1\LL B_L)$ is {\em $\lambda$-top-concave} if it has the following very strong property:
for every $Y\subseteq B_1\cup\cdots\cup B_L$, there do not exist $h_1,h_2,h_3$ with $1\le h_1<h_2<h_3 \le K$,
such that $Y$ $\lambda$-covers $A_{h_2}$ and $Y$ $\lambda$-misses $A_{h_1}$ and $A_{h_3}$. We define
{\em $\lambda$-bottom-concave} similarly. We say $\mathcal{P}$ is {\em $\lambda$-concave} if it is both $\lambda$-top-concave
and $\lambda$-bottom-concave.

A bigraph $G$ is {\em balanced} if $|V_1(G)|=|V_2(G)|$; and
we say a parade $(A_1\LL A_K; B_1\LL B_L)$ is {\em balanced} if $K=L$ and all its blocks have the same cardinality (that is, 
$|A_1|=|B_1|$).

The following theorem, which is the main step in the proof, says that we can find a rooted
tree bigraph $T$ in any parade that is sufficiently well-behaved.  In later sections, we will show that it is possible to 
find such a parade.

\begin{thm}\label{usemonotone}
Let $\delta\ge 2$ and $\eta\ge 0$ be integers.
Let $T$ be a rooted tree bigraph, such that every vertex has degree at most $\delta+1$, the root has
degree at most $\delta$,
and every path
from root to leaf has length less than $\eta$. Let $\tau=\delta^{\eta+1}$, and
let
$0<\lambda \le 2^{-30\delta}\delta^{-1-\eta}$.
Let $G$ be a balanced bigraph with a balanced parade
$\mathcal{P}$ of length $(K,K)$ where $K=(32\delta+4)\tau+2$, such that
$\mathcal{P}$ is $\lambda$-concave, $(2^{-30\delta},\tau)$-support-invariant and $\tau$-support-uniform.
Let $\mathcal{P}$ have width $(W,W)$.
If $G$ is $\vare$-coherent where $\vare\le 2^{-30\delta}$, then there is a $\mathcal{P}$-rainbow copy of $T$.
\end{thm}
\Proof
We will be looking at rooted tree bigraphs in which every vertex has degree exactly $\delta+1$ except the root and the leaves; but the root
might have degree different from $\delta$, and not all paths from root to leaf will necessarily have the same length. Let us first set
up some notation for such trees.

If $a\ge 2$ is an integer, let $T(a,0)$ be the rooted tree bigraph $H$ with $|V_1(H)|=1$ and $V_2(H)=\emptyset$
(thus, $a$ is irrelevant, but this will be convenient).
If $a\ge 2$ and $b\ge 1$ are integers, let $T(a,b)$ be the rooted tree bigraph $H$ with root in $V_1(H)$, such that the root has degree $a$, every vertex different
from the root has degree $a+1$ or $1$, and every path from root to leaf has length exactly $b$. We denote
the transpose of $T(a,b)$
by $\tilde{T}(a,b)$.
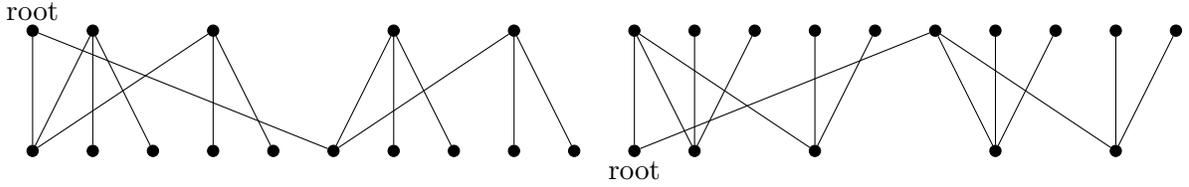
\begin{figure}[H]
\centering

\begin{tikzpicture}[scale=0.8,auto=left]
\tikzstyle{every node}=[inner sep=1.5pt, fill=black,circle,draw]
\node (a) at (1,2) {};
\node (b1) at (1,0) {};
\node (b2) at (6,0) {};
\node (c1) at (2,2) {};
\node (c2) at (4,2) {};
\node (c3) at (7,2) {};
\node (c4) at (9,2) {};
\node (d1) at (2,0) {};
\node (d2) at (3,0) {};
\node (d3) at (4,0) {};
\node (d4) at (5,0) {};
\node (d5) at (7,0) {};
\node (d6) at (8,0) {};
\node (d7) at (9,0) {};
\node (d8) at (10,0) {};

\foreach \from/\to in {a/b1,a/b2,b1/c1,b1/c2,b2/c3,b2/c4,c1/d1,c1/d2,c2/d3,c2/d4,c3/d5,c3/d6,c4/d7,c4/d8}
\draw [-] (\from) -- (\to);

\tikzstyle{every node}=[]
\draw[above] (a) node []           {root};

\begin{scope}[shift ={(10,2)}]
\tikzstyle{every node}=[inner sep=1.5pt, fill=black,circle,draw]
\node (a) at (1,-2) {};
\node (b1) at (1,0) {};
\node (b2) at (6,0) {};
\node (c1) at (2,-2) {};
\node (c2) at (4,-2) {};
\node (c3) at (7,-2) {};
\node (c4) at (9,-2) {};
\node (d1) at (2,0) {};
\node (d2) at (3,0) {};
\node (d3) at (4,0) {};
\node (d4) at (5,0) {};
\node (d5) at (7,0) {};
\node (d6) at (8,0) {};
\node (d7) at (9,0) {};
\node (d8) at (10,0) {};

\foreach \from/\to in {a/b1,a/b2,b1/c1,b1/c2,b2/c3,b2/c4,c1/d1,c1/d2,c2/d3,c2/d4,c3/d5,c3/d6,c4/d7,c4/d8}
\draw [-] (\from) -- (\to);

\tikzstyle{every node}=[]
\draw[below] (a) node []           {root};

\end{scope}
\end{tikzpicture}
\caption{$T(2,3)$ and $\tilde{T}(2,3)$.} \label{fig:T(2,3)}
\end{figure}

Now let $a_1,a_2, b_1,b_2\ge 0$ be integers. Let $T(a_1,b_1,a_2,b_2)$ be the rooted tree bigraph obtained from
$a_1$ copies of $\tilde{T}(\delta,b_1)$ and $a_2$ copies of $\tilde{T}(\delta,b_2)$, all pairwise disjoint, by adding a new root
adjacent to all the old roots. We denote its transpose by $\tilde{T}(a_1,b_1,a_2,b_2)$.

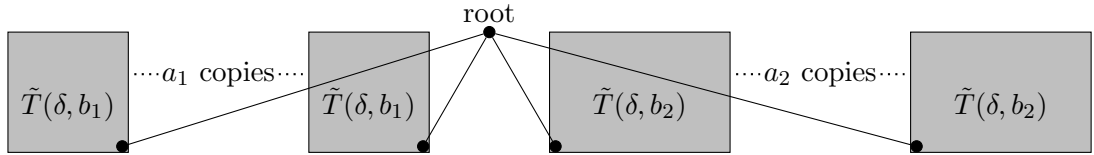
\begin{figure}[h!]
\centering

\begin{tikzpicture}[scale=0.8,auto=left]

\draw[fill= lightgray] (-1,0) rectangle (-3,2);
\draw[fill= lightgray] (-6,0) rectangle (-8,2);
\draw[fill= lightgray] (1,0) rectangle (4,2);
\draw[fill= lightgray] (7,0) rectangle (10,2);

\tikzstyle{every node}=[inner sep=1.5pt, fill=black,circle,draw]
\node (a) at (0,2) {};
\node (b) at (-1.1,0.1) {};
\node (c) at (-6.1,0.1) {};

\node (d) at (1.1,0.1) {};
\node (e) at (7.1,0.1) {};

\draw[-] (a)--(b);
\draw[-] (a)--(c);
\draw[-] (a)--(d);
\draw[-] (a)--(e);

\tikzstyle{every node}=[]
\draw[above] (a) node []           {root};
\node at (-2,.8) {$\tilde{T}(\delta,b_1)$};
\node at (-7,.8) {$\tilde{T}(\delta,b_1)$};
\node at (-4.5,1.3) {$a_1$ copies};

\node at (2.5,.8) {$\tilde{T}(\delta,b_2)$};
\node at (8.5,.8) {$\tilde{T}(\delta,b_2)$};
\node at (5.5,1.3) {$a_2$ copies};

\draw[dotted,thick] (-5.9,1.3)--(-5.5,1.3);
\draw[dotted,thick] (-3.5,1.3)--(-3.1,1.3);
\draw[dotted,thick] (4.1,1.3)--(4.5,1.3);
\draw[dotted,thick] (6.5,1.3)--(6.9,1.3);

\end{tikzpicture}
\caption{$T(a_1,b_1,a_2,b_2)$.} \label{T(a_1,b_1,a_2,b_2)}
\end{figure}
Let $\mathcal{P}=(A_1\LL A_K; B_1\LL B_K)$ as in the theorem, where $K=(32\delta+4)\tau+2$, and let its width be $(W,W)$.
We observe that, since $G$ is $\vare$-coherent, fewer than $\vare W$ vertices in $A_1$ have no neighbour in $B_{i}$ for 
$1\le i  \le K$; and since $\delta\le K$, and $\vare\delta <1$, there is a vertex in $A_1$ with a neighbour in each of
$B_{1}\LL B_{\delta}$. Hence there is a copy of $T(\delta,1)$ that is $\mathcal{P}$-left-rainbow and $\mathcal{P}$-right-rainbow;
and similarly there is a copy of $\tilde{T}(\delta,1)$ that is $\mathcal{P}$-left-rainbow and $\mathcal{P}$-right-rainbow.
Consequently we may choose $\rho,\tilde{\rho},\sigma,\tilde{\sigma}$ as below.
\begin{itemize}
\item Choose $\rho\ge 1$ maximum such that there is a copy of $T(\delta,\rho)$ that is $\mathcal{P}$-left-rainbow (and therefore there is a $\mathcal{P}$-left-rainbow copy of $T(0,\rho,\delta,\rho-1)$);
and choose $\phi\ge 0$ maximum such that there is a copy of $T(\phi,\rho,\delta-\phi,\rho-1)$ that is $\mathcal{P}$-left-rainbow.
\item Choose $\tilde{\rho}\ge 1$ maximum such that there is a copy of $\tilde{T}(\delta,\tilde{\rho})$ that is
$\mathcal{P}$-left-rainbow; and choose $\tilde{\phi}$ maximum such that there is a copy of
$\tilde{T}(\tilde{\phi},\tilde{\rho},\delta-\tilde{\phi},\tilde{\rho}-1)$ that is $\mathcal{P}$-left-rainbow.
\item Choose $\sigma\ge 1$ maximum such that there is a copy of $T(\delta,\sigma)$ that is $\mathcal{P}$-right-rainbow;
and choose $\psi\ge 0$ maximum such that there is a copy of $T(\psi,\sigma,\delta-\psi,\sigma-1)$ that is $\mathcal{P}$-right-rainbow.
\item Choose $\tilde{\sigma}\ge 1$ maximum such that there is a copy of $\tilde{T}(\delta,\tilde{\sigma})$ that is
$\mathcal{P}$-right-rainbow; and choose $\tilde{\psi}$ maximum such that there is a copy of
$\tilde{T}(\tilde{\psi},\tilde{\sigma},\delta-\tilde{\psi},\tilde{\sigma}-1)$ that is $\mathcal{P}$-right-rainbow.
\end{itemize}
We suppose for a contradiction that there is no $\mathcal{P}$-rainbow copy of $T(\delta,\eta)$ or its transpose, and so
$\rho,\tilde{\rho},\sigma,\tilde{\sigma}<\eta$. Also $\phi<\delta$, since otherwise $T(\phi,\rho,\delta-\phi,\rho-1)$
contains $T(\delta,\rho+1)$, contrary to the maximality of $\rho$; and similarly $\tilde{\phi},\psi,\tilde{\psi}<\delta$.

Let us partition $\{2\LL K-1\}$ into $16\delta+2$ intervals,
each of length $2\tau$, and numbered in order. Thus, $I_0\LL I_{16\delta+1}$ are pairwise disjoint subsets of
$\{2\LL K-1\}$, with union $\{2\LL K-1\}$, each of cardinality
$2\tau$, and such that $x<y$ for all $i,j\in \{0\LL 16\delta+1\}$ with $i<j$ and all $x\in I_i$ and $y\in I_j$.

Let $a,b,c,d,\tilde{a}, \tilde{b},\tilde{c},\tilde{d}\ge 0$ be integers, all at most $\delta$, and
let $e$ be the sum of these eight integers.
Let $I$ be the union of $\{1,K\}$ and the sets $I_j$ for all $j\in \{0\LL 8\delta-e\}\cup \{8\delta+e+1\LL 16\delta+1\}$.
Let $\mathcal{P}'=(A_i'\;(i\in I);B_i'\;(i\in I))$ be a balanced minor of $\mathcal{P}$, where $A_i'\subseteq A_i$ and
$B_i'\subseteq B_i$ for $i\in I$.
We say that $\mathcal{P}'$ is {\em $(a,b,c,d,\tilde{a}, \tilde{b},\tilde{c},\tilde{d})$-anchored}
if
the following hold:
\begin{itemize}
\item $\mathcal{P}'$ has width at least $2^{-3e}(W,W)$.
\item $Y\subseteq \bigcup_{i\in \{1\LL K\}\setminus I}A_i\cup B_i$, and $Y$ is anticomplete to $A_{i}'\cup B_{i}'$ for all $i\in I\setminus \{1,K\}$.
\item For every $v\in A_{1}'$ there is a $\mathcal{P}$-left-rainbow copy of $T(\tilde{a},\tilde{\rho}, \tilde{b}, \tilde{\sigma})$ in $G[Y\cup \{v\}]$ with root $v$.
\item For every $v\in B_{1}'$ there is a $\mathcal{P}$-left-rainbow copy of $\tilde{T}(a,\rho,b,\sigma)$ in $G[Y\cup \{v\}]$ with root $v$.
\item For every $v\in A_{K}'$ there is a $\mathcal{P}$-right-rainbow copy of $T(\tilde{c},\tilde{\rho}, \tilde{d}, \tilde{\sigma})$ in $G[Y\cup \{v\}]$ with root $v$.
\item For every $v\in B_{K}'$ there is a $\mathcal{P}$-right-rainbow copy of $\tilde{T}(c,\rho,d,\sigma)$ in $G[Y\cup \{v\}]$ with root $v$.
\end{itemize}
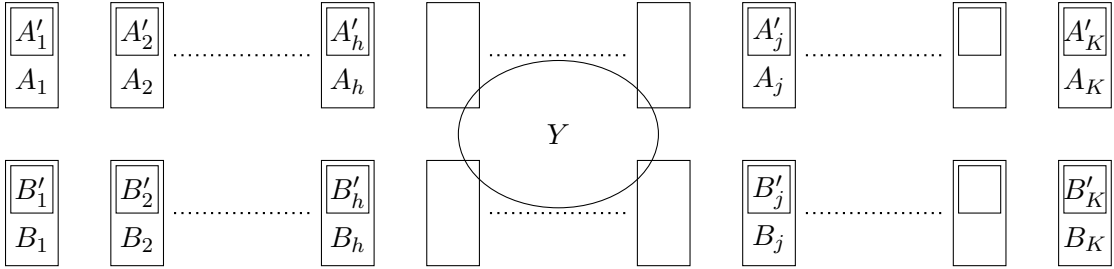
\begin{figure}[h!]
\centering

\begin{tikzpicture}[scale=0.7,auto=left]
\tikzstyle{every node}=[]

\draw (0,1) rectangle (1,3);
\draw (2,1) rectangle (3,3);
\draw (6,1) rectangle (7,3);
\draw (8,1) rectangle (9,3);
\draw (12,1) rectangle (13,3);
\draw (14,1) rectangle (15,3);
\draw (18,1) rectangle (19,3);
\draw (20,1) rectangle (21,3);

\draw (0.1,2) rectangle (.9,2.9);
\draw (2.1,2) rectangle (2.9,2.9);
\draw (6.1,2) rectangle (6.9,2.9);
\draw (14.1,2) rectangle (14.9,2.9);
\draw (18.1,2) rectangle (18.9,2.9);
\draw (20.1,2) rectangle (20.9,2.9);

\draw[dotted, thick] (3.2,2) -- (5.8,2);
\draw[dotted, thick] (9.2,2) -- (11.8,2);
\draw[dotted, thick] (15.2,2) -- (17.8,2);

\node at (.5,1.5) {$A_1$};
\node at (2.5,1.5) {$A_2$};
\node at (6.5,1.5) {$A_{h}$};
\node at (14.5,1.5) {$A_{j}$};
\node at (20.5,1.5) {$A_K$};

\node at (.5,2.4) {$A_1'$};
\node at (2.5,2.4) {$A_2'$};
\node at (6.5,2.4) {$A_h'$};
\node at (14.5,2.4) {$A_j'$};
\node at (20.5,2.4) {$A_K'$};

\begin{scope}[shift={(0,-3)}]
\draw (0,1) rectangle (1,3);
\draw (2,1) rectangle (3,3);
\draw (6,1) rectangle (7,3);
\draw (8,1) rectangle (9,3);
\draw (12,1) rectangle (13,3);
\draw (14,1) rectangle (15,3);
\draw (18,1) rectangle (19,3);
\draw (20,1) rectangle (21,3);

\draw (0.1,2) rectangle (.9,2.9);
\draw (2.1,2) rectangle (2.9,2.9);
\draw (6.1,2) rectangle (6.9,2.9);
\draw (14.1,2) rectangle (14.9,2.9);
\draw (18.1,2) rectangle (18.9,2.9);
\draw (20.1,2) rectangle (20.9,2.9);

\draw[dotted, thick] (3.2,2) -- (5.8,2);
\draw[dotted, thick] (9.2,2) -- (11.8,2);
\draw[dotted, thick] (15.2,2) -- (17.8,2);

\node at (.5,1.5) {$B_1$};
\node at (2.5,1.5) {$B_2$};
\node at (6.5,1.5) {$B_{h}$};
\node at (14.5,1.5) {$B_{j}$};
\node at (20.5,1.5) {$B_K$};

\node at (.5,2.4) {$B_1'$};
\node at (2.5,2.4) {$B_2'$};
\node at (6.5,2.4) {$B_h'$};
\node at (14.5,2.4) {$B_j'$};
\node at (20.5,2.4) {$B_K'$};
\end{scope}
\node at (10.5,.5) {$Y$};
\draw (10.5,.5) ellipse (1.9 and 1.4);

\end{tikzpicture}
\caption{Figure for ``anchored''. Here, $h=2\tau(8\delta-e)+1$, and $j=2\tau(8\delta+e)+2$.} \label{fig:anchored}
\end{figure}

Let $|V_1(G)|=|V_2(G)|=n$. Choose $a,b,c,d,\tilde{a}, \tilde{b},\tilde{c},\tilde{d}$ with maximum sum
such that some balanced minor is $(a,b,c,d,\tilde{a}, \tilde{b},\tilde{c},\tilde{d})$-anchored (this is possible since
$\mathcal{P}$ is $(0\LL 0)$-anchored),
and let $I, \mathcal{P}', a,b,c,d,\tilde{a}, \tilde{b},\tilde{c},\tilde{d}, e$ and $Y$ be as above.
Let $\mathcal{P}'$ have width $(W', W')$. Since $e\le 8\delta$, it follows that
$W'\ge 2^{-24\delta}W$.

For a bigraph $H$, we call the pair $(|V_1(H)|,|V_2(H)|)$ its {\em part size}.
Let 
$$T(\phi,\rho,\delta-\phi,\rho-1), 
\tilde{T}(\tilde{\phi},\tilde{\rho},\delta-\tilde{\phi},\tilde{\rho}-1),
T(\psi,\sigma,\delta-\psi,\sigma-1), 
\tilde{T}(\tilde{\psi},\tilde{\sigma},\delta-\tilde{\psi},\tilde{\sigma}-1)$$
have part size $(p_1,q_1), (p_2,q_2), (p_3,q_3), (p_4,q_4)$ respectively. Each of these eight numbers is at most $\tau$.
Choose
\begin{itemize}
\item disjoint subsets $P_1,P_2$ of $I_{8\delta-e}$ of cardinality $p_1,p_2$ respectively, with the smallest member of 
$I_{8\delta-e}$ in $P_1$; 
\item disjoint subsets $P_3, P_4$ of $I_{8\delta+e+1}$ of cardinality $p_3, p_4$ respectively, with the 
largest member of $I_{8\delta+e+1}$ in $P_3$;
\item disjoint subsets $Q_1,Q_2$ of $I_{8\delta-e}$ of cardinality $q_1,q_2$ respectively,
with the smallest member of $I_{8\delta-e}$ in $Q_2$; 
\item disjoint subsets $Q_3, Q_4$ of $I_{8\delta+e+1}$ of cardinality $q_3, q_4$ respectively,
with the largest member of $I_{8\delta+e+1}$ in $Q_4$.
\end{itemize}
Let $r=\lceil (2^{-24\delta}-2^{-30\delta}) W \rceil$. It follows that
$r\ge (2^{6\delta}-1)2^{-30\delta}W \ge (2^{12}-1)2^{-30\delta}W$ since $\delta\ge 2$.
\\
\\
(1) {\em There are $r$ copies $L_1\LL L_r$ of $T(\phi,\rho,\delta-\phi,\rho-1)$,
pairwise vertex-disjoint and each $\mathcal{P}'$-left-rainbow with support $(P_1,Q_1)$.}
\\
\\
Since there is a copy of $T(\phi,\rho,\delta-\phi,\rho-1)$ that is $\mathcal{P}$-left-rainbow, and $\mathcal{P}$ is
$\tau$-support-uniform,
there is such a copy that is $\mathcal{P}$-left-rainbow with support $(P_1,Q_1)$.
Choose $r'\le r$ maximum such that there are $r'$ copies of $T(\phi,\rho,\delta-\phi,\rho-1)$,
pairwise vertex-disjoint and each $\mathcal{P}'$-left-rainbow with support $(P_1,Q_1)$. Suppose that $r'<r$. By
removing the vertices of these copies from the blocks of $\mathcal{P}'$ that contain them, and removing $r'$ vertices arbitrarily
from every other block of $\mathcal{P}'$,
we obtain a minor of $\mathcal{P}'$, and hence of $\mathcal{P}$, with width $(W'-r', W'-r')$,
relative to which there is no left-rainbow copy of $T(\phi,\rho,\delta-\phi,\rho-1)$ with support $(P_1,Q_1)$.
But $\mathcal{P}'$ is
$(2^{-30\delta},\tau)$-support-invariant,
and so $W'-r'<2^{-30\delta}W$. Since $W' \ge 2^{-24\delta}W$, it follows that
$r'>(2^{-24\delta}-2^{-30\delta}) W $, contradicting that $r'<r$. This proves (1).

\bigskip

Similarly, we deduce:
\begin{itemize}
\item There are $r$ copies $\tilde{L}_1\LL \tilde{L}_r$ of $\tilde{T}(\tilde{\phi},\tilde{\rho},\delta-\tilde{\phi},\tilde{\rho}-1)$,
pairwise vertex-disjoint and each $\mathcal{P}'$-left-rainbow with support $(P_2,Q_2)$.
\item There are $r$ copies $M_1\LL M_r$ of $T(\psi,\sigma,\delta-\psi,\sigma-1)$,
pairwise vertex-disjoint and each $\mathcal{P}$-right-rainbow with support $(P_3,Q_3)$.
\item There are $r$ copies $\tilde{M}_1\LL \tilde{M}_r$ of  $\tilde{T}(\tilde{\psi},\tilde{\sigma},\delta-\tilde{\psi},\tilde{\sigma}-1)$,
pairwise vertex-disjoint and each $\mathcal{P}$-right-rainbow relative to $\mathcal{P}'$ with support $(P_4,Q_4)$.
\end{itemize}

For $1\le i\le r$, let $H_i$ be the disjoint union of $L_i, \tilde{L}_i, M_i$ and $\tilde{M}_i$.
Thus $H_i$ is $\mathcal{P}'$-rainbow.  If $h$ is the smallest member of $I_{8\delta-e}$, then the root of $L_i$ belongs to $A_{h}'$ and the root
of $\tilde{L}_i$ belongs to $B_{h}'$; and if $h'$ is the largest member of $I_{8\delta+e+1}$, then the root of $M_i$
belongs to $A_{h'}'$, and the root of $\tilde{M}_i$ belongs to $B_{h'}'$. We call these four roots the {\em extremes} of $H_i$.
Let $v\in A_1'\cup B_1'\cup A_K'\cup B_K'$, and $1\le i\le r$. If $v\in A_1'\cup A_K'$ then every vertex of $H_i$ adjacent to
$v$ belongs to $V_2(H_i)$, and if $v\in B_1'\cup B_K'$ then every vertex of $H_i$ adjacent to $v$ belongs to $V_1(H_i)$.
In particular,
there are at most two extremes of $H_i$ that are adjacent to $v$. We say
\begin{itemize}
\item  $v$ {\em meets} $H_i$ if $v$ is adjacent to some vertex of
$H_i$;
\item $v$ meets $H_i$ {\em internally} if $v$ is adjacent to some vertex of $H_i$
that is not an extreme (and possibly $v$ is also adjacent to one or two extremes);
\item $v$ meets $H_i$ {\em properly} if $v$ is adjacent to one or two extremes, and to no other
vertices of $H_i$, that is, if $v$ meets $H_i$ and does not meet $H_i$ internally.
\end{itemize}

For $X\subseteq A_1'\cup B_1'\cup A_K'\cup B_K'$, let $f(X)$ be the number of $i\in \{1\LL r\}$ such that some vertex in $X$ meets $H_i$,
and let $g(X)$ be the number of $i\in \{1\LL r\}$ such that some vertex in $X$ meets $H_i$ internally. Choose
$D\subseteq A_1'\cup B_1'\cup A_K'\cup B_K'$
maximal such that $f(D)\le r/2$ and $g(D)\ge f(D)/8$. For $i = 1,2$, let $D_i=D\cap V_i(G)$.
\\
\\
(2) {\em $D_1,D_2$ both have cardinality less than $\vare n$; and
$f(D)\le r/2-\vare n$.}
\\
\\
Let $h$ be the smallest element of $I_{8\delta-e}$.
There are at least $r/2$ vertices in $A_{h}$
with no neighbour in $D_2$ (the roots of the trees $L_i$ such that no vertex in $D$ meets $H_i$). Since
$r/2\ge \vare n$ and $G$ is $\vare$-coherent, it follows that
$|D_2|< \vare n$; and since $r/2\ge \lambda W$,
we deduce that $D_2$ $\lambda$-misses $A_{h}$.
Similarly it $\lambda$-misses $A_{h'}$, where $h'$ is the largest member of $I_{8\delta+e+1}$; and since
$\mathcal{P}$ is $\lambda$-concave, $D_2$ does not $\lambda$-cover any of the sets $A_{j} (h<j<h')$.
In particular, since $|I_{8\delta-e}|=|I_{8\delta+e+1}|=2\tau$, there are at most $\lambda (4\tau)W$ vertices in
$$\bigcup \left(A_j:j\in (I_{8\delta-e}\setminus \{h\}) \cup (I_{8\delta+e+1}\setminus \{h'\})\right)$$
that have neighbours in $D_2$. It follows that $g(D_2)\le 4\lambda \tau W$. The same holds
for $D_1$; and since
$g(D)\le g(D_1)+g(D_2),$
it follows that $g(D)\le 8\lambda \tau W$. But $f(D)\le 8g(D)$, and so
$$f(D)\le 64\lambda\tau W 
\le 64(2^{-30\delta}\delta^{-1-\eta})(\delta^{\eta+1})((2^{12}-1)^{-1}2^{30\delta}r)=
64(2^{12}-1)^{-1}r\le r/4\le r/2-\vare n$$
since $\lambda \le 2^{-30\delta}\delta^{-1-\eta}$, and
$r\ge (2^{12}-1)2^{-30\delta} W$, and $r\ge 4\vare n$.
This proves (2).
\bigskip

Let $F$ be the set of vertices in $(A_1'\cup A_K'\cup B_1'\cup B_K')\setminus D$ that meet one of $H_1\LL H_r$.
\\
\\
(3) {\em At most $2\vare n$ vertices in $A_1'\cup A_K'$ do not belong to $F$, and the same for $B_1'\cup B_K'$.}
\\
\\
Since $r\ge \vare n$,
and $G$ is $\vare$-coherent, there are fewer than $\vare n$ vertices in $A_1'\cup A_K'$ that have no neighbour in any of
$H_1\LL H_r$. All the other vertices in $A_1'\cup A_K'$ belong to either $D_1$ or $F$, and only at most $\vare n$ belong to $D_1$, by (2).
Consequently at most $2\vare n$ vertices in $A_1'\cup A_K'$ do not belong to $F$, and the same for $B_1'\cup B_K'$.
This proves (3).

\bigskip

Let $C$ be the set of all $i\in \{1\LL r\}$ such that $D$ is anticomplete to $V(H_i)$. Thus $|C|=r-f(D)$.
\\
\\
(4) {\em For each $v\in F$, the number of $i\in C$ such that $v$ meets $H_i$ internally is less than $1/8$
of the number of $i\in C$ such that $v$ meets $H_i$.}
\\
\\
Since $f(D)\le r/2-\vare n$ by (2), it follows that
$f(D\cup \{v\})\le r/2$, and the maximality of $D$ implies that
$g(D\cup\{v\})< f(D\cup \{v\})/8$. Since $g(D)\ge f(D)/8$, this proves (4).

\bigskip
Let us say $v\in F$ is {\em happy} if  $v$ meets $H_i$ properly, where $i\in C$ is minimum such that $v$ meets $H_i$.
\\
\\
(5) {\em We may assume that at least half of the vertices in
each of the sets $A_1'\cap F, A_K'\cap F, B_1'\cap F, B_K'\cap F$ are happy.}
\\
\\
Let $v\in F$, and take a linear order of $C$.
If we choose the linear order uniformly at random,
the probability that $v$ is happy is more than $7/8$, by (4); and so the expected number of happy vertices in $A_1'$
is more than $7/8 |A_1'|$.
Hence the probability that at least half the vertices in $A_1'$ are happy is more than $3/4$
(because if it were at most $3/4$ then the expected number of happy vertices would be at most $(3/4)|A_1'|+(1/4)|A_1'|/2$, which
is too small). Similarly the probability that at least half the vertices in $A_K'$ are happy is more than $3/4$, and the same
for $B_1', B_K'$; and so there is positive probability that all four events happen, that is, for some linear order of $C$,
at least half the vertices of each of $A_1', A_K', B_1', B_K'$ are happy. Renumber $C$ in this order; then
(5) holds.

\bigskip
Let $X$ be the set of all happy vertices in $F$. For each $v\in X$, choose $i\in C$ minimum such that $v$ meets $H_i$.
We call $i$ the {\em happiness} of $v$.
Now $|F\cap A_1'|\ge W-2\vare n$ by (3), and so
$|X\cap A_1'|\ge W'/2-\vare n\ge 3W'/8$.
Thus we may
choose $m\le |C|$ minimum such that one of $A_1'\cap X, A_K'\cap X, B_1'\cap X,B_K'\cap X$ contains at least $W'/4$
vertices with happiness at most $m$.
Let $Y'=V(H_1)\cup\cdots\cup V(H_m)$ and $I'=I\setminus (I_{8\delta-e}\cup I_{8\delta+e+1})$.
\\
\\
(6) {\em For each $i\in I'\setminus \{1,K\}$, there is a subset $A_i''$ of $A_i'$, and a subset $B_i''$ of $B_i'$, both anticomplete to $Y'$,
and both of cardinality $\lceil W'/8\rceil$.}
\\
\\
Let $i\in I'\setminus \{1,K\}$; we will show that $A_i'$ has a subset with the desired properties.
Let $j\in I_{8\delta-e}\cup I_{8\delta+e+1}$.
From the choice of $m$, fewer than $W'/4$ vertices in $X\cap A_1'$ have happiness less
than $m$; and so at most $W'/4 +2\vare n$
have happiness at most $m$, since those with happiness exactly $m$ are adjacent to one of the roots of $\tilde{L}_m, \tilde{M}_m$.
Since $|X\cap A_1'|\ge 3W'/8$, there are at least $W'/8-2\vare n$
vertices in $X\cap A_1'$ that have no neighbour in $Y'$, and in particular have no neighbour in $B_j'\cap Y'$. Since
$W'/8-2\vare n>\lambda W$,
$B_j'\cap Y'$
$\lambda$-misses $A_1$. By the same argument $B_j'\cap Y'$ $\lambda$-misses $A_K$, and so does not $\lambda$-cover
$A_i$, since $\mathcal{P}$ is $\lambda$-concave.
Consequently, for each $j\in I_{8\delta-e}\cup I_{8\delta+e+1}$,  there are at most $\lambda W$ vertices in $A_i'$ with a neighbour in
$B_j'\cap Y'$; and so there are at most $4\tau\lambda W$ vertices in $A_i'$ with a neighbour in $Y'$.
Since $|A_i'|= W'$ and $W'-4\tau\lambda W n\ge W'/8$, this proves that $A_i'$ has a subset with the
desired properties.
From the symmetry under taking transpose, this proves (6).

\bigskip
We chose $m$ such that at least one of $A_1'\cap X, A_K'\cap X, B_1'\cap X,B_K'\cap X$ contains at least $W'/4$
vertices with happiness at most $m$; and from the symmetry, we may assume that $B_1'\cap X$ contains
at least $W'/4$
vertices with happiness at most $m$. If $v\in B_1'\cap X$ has happiness at most $m$, let $i$ be its happiness; then
$v$ is adjacent to either the root of $L_i$ or the root of $M_i$.
Choose $Z\subseteq B_1'\cap X$ with $|Z|\ge W'/8$ such that either every vertex $v\in Z$ is adjacent to the root of
$L_i$, where $i$ is the happiness of $v$, or every $v\in Z$ is adjacent to the root of
$M_i$, where $i$ is the happiness of $v$.
Choose $B_1''\subseteq B_1'$ of cardinality  $\lceil W'/8\rceil$ with $B_1''\subseteq Z$ (this is possible since
$|Z|\ge W'/8$). Choose $A_1''\subseteq A_1'$, $A_K''\subseteq A_K'$, and $B_K''\subseteq B_K'$, all of cardinality
$\lceil W'/8\rceil$.
Let $\mathcal{P}''=(A_i''\;(i\in I'');B_i''\;(i\in I'))$. Then $\mathcal{P}''$ is a balanced minor of $\mathcal{P}$.
Its width is at least $W'/8$.

Suppose first that every vertex $v\in Z$ is adjacent to the root of
$L_i$, where $i$ is the happiness of $v$.
Choose $v\in Z$, and let $i$ be its happiness. Let $u$ be the root of $L_i$.
Since
$\mathcal{P}'$ is $(a,b,c,d,\tilde{a}, \tilde{b},\tilde{c},\tilde{d})$-anchored, there is a copy of $\tilde{T}(a,\rho, b, \sigma)$
in $G[Y\cup \{v\}]$ with root $v$, $\mathcal{P}$-left-rainbow, say $S$.

Suppose that $a=\delta$.  Then $S$ has a rooted subtree $S'$ with root $v$, isomorphic
to $\tilde{T}(\delta,\rho)$. But
$L_i$ is a copy of $T(\phi,\rho,\delta-\phi,\rho-1)$; and so the union of $S'$, $L_i$, and the edge $uv$, with root $u$, 
is a $\mathcal{P}$-left-rainbow copy of $T(\phi+1,\rho,\delta-\phi,\rho-1)$.
From the choice of $\rho$, $\phi+1<\delta$;
and so there is a $\mathcal{P}$-left-rainbow copy of $T(\phi+1,\rho,\delta-\phi-1,\rho-1)$, contrary to the choice of $\phi$.
This proves that $a<\delta$.

By taking the union of $S$ and an appropriate subtree of $L_i$
and the edge $uv$, we obtain a copy of $T(a+1,\rho, b, \sigma)$ in $G[Y\cup Y'\cup \{v\}]$, with root $v$. This holds for each
$v\in Z$.
We claim that $\mathcal{P}''$ 
is $(a+1,b,c,d,\tilde{a}, \tilde{b},\tilde{c},\tilde{d})$-anchored. Since its width is at least $2^{-3e}(W,W)$, it suffices to check that:
\begin{itemize}
\item $Y\cup Y'\subseteq \bigcup\left(\bigcup_{i\in I_j}A_i\cup B_i:8\delta-(e+1)+1\le j\le 8\delta+(e+1)\right)$, and $Y\cup Y'$
is anticomplete to $A_{i}''\cup B_{i}''$ for all $i\in I\setminus \{1,K\}$.
\item For every $v\in A_{1}''$ there is a copy of $T(\tilde{a},\tilde{\rho}, \tilde{b}, \tilde{\sigma})$ in $G[Y\ \cup Y'\cup \{v\}]$ with root $v$, $\mathcal{P}$-left-rainbow.
\item For every $v\in B_{1}''$ there is a copy of $\tilde{T}(a+1,\rho,b,\sigma)$ in $G[Y\cup Y'\cup \{v\}]$ with root $v$, $\mathcal{P}$-left-rainbow.
\item For every $v\in A_{K}''$ there is a copy of $T(\tilde{c},\tilde{\rho}, \tilde{d}, \tilde{\sigma})$ in $G[Y\cup Y'\cup \{v\}]$ with root $v$, $\mathcal{P}$-right-rainbow.
\item For every $v\in B_{K}''$ there is a copy of $\tilde{T}(c,\rho,d,\sigma)$ in $G[Y\cup Y'\cup \{v\}]$ with root $v$, $\mathcal{P}$-right-rainbow.
\end{itemize}
The first of these holds from the choice of $Y'$ and of the sets $A_i'', B_i''$.
We have just seen that the third holds;
and the other three statements are true because $\mathcal{P}'$ is
$(a,b,c,d,\tilde{a}, \tilde{b},\tilde{c},\tilde{d})$-anchored. But this contradicts the maximality of $e$.

This completes the case when every vertex $v\in Z$ is adjacent to the root of
$L_i$, where $i$ is the happiness of $v$; so now we may assume that every vertex $v\in Z$ is adjacent to the root of
$M_i$, where $i$ is the happiness of $v$. In fact this second case is the same as the first case, as can be seen by
reversing the numbering of $A_1\LL A_K$; but checking that this symmetry argument is valid seems more difficult than repeating the argument
for the first case, so we will just repeat the argument for the first case.

Choose $v\in Z$, and let $i$ be its happiness. Since
$\mathcal{P}'$ is $(a,b,c,d,\tilde{a}, \tilde{b},\tilde{c},\tilde{d})$-anchored, there is a $\mathcal{P}$-left-rainbow 
copy of $\tilde{T}(a,\rho, b, \sigma)$
in $G[Y\cup \{v\}]$ with root $v$, say $S$. Let $u$ be the root of $M_i$.

Suppose that $b=\delta$.
Then $S$ has a rooted subtree $S'$ with root $v$ isomorphic
to $\tilde{T}(\delta,\sigma)$.  But
$M_i$ is a copy of $T(\psi,\sigma,\delta-\psi,\sigma-1)$; and so the union of $S'$, $M_i$, and the edge $uv$, 
with root $u$, is a $\mathcal{P}$-left-rainbow copy of $T(\psi+1,\sigma,\delta-\psi,\sigma-1)$.
From the choice of $\sigma$, $\psi+1<\delta$;
and so there is a $\mathcal{P}$-left-rainbow copy of $T(\psi+1,\sigma,\delta-\psi-1,\sigma-1)$, contrary to the choice of $\psi$.
This proves that $b<\delta$.

By taking the union of $S$ with an appropriate subtree of $M_i$
and the edge $uv$, we obtain a copy of $T(a,\rho, b+1, \sigma)$ in $G[Y\cup Y'\cup \{v\}]$.
We claim that $\mathcal{P}''$ is
is $(a,b+1,c,d,\tilde{a}, \tilde{b},\tilde{c},\tilde{d})$-anchored, and the argument is as before.
But this contradicts the maximality of $e$.

This contradiction shows that there is a $\mathcal{P}$-rainbow copy of $T(\delta,\eta)$ or its transpose, and since both these
contain copies of $T$, possibly with different root (because all paths in $T$ from root to leaf have length less than $\eta$),
there is also
a $\mathcal{P}$-rainbow copy of $T$.
This proves \ref{usemonotone}.~\bbox

\section{Unbalanced parades}

We have finished the difficult part of the paper; now we just have to apply \ref{usemonotone}. One problem is that \ref{usemonotone}
applies only to balanced parades in balanced bigraphs, and it would be easier to use without that restriction.
In this section we deduce a version of \ref{usemonotone} without the balancedness restrictions.

For clarity, in what follows we say ``$G$-adjacent'' to mean adjacent in $G$, and define ``$G$-neighbour'', ``$G$-anticomplete'' and so on,
similarly.
Let $G$ be a bigraph, and let $a,b>0$ be integers. For each $u\in V_1(G)$ take a set $M_u$ of $a$ new vertices,
and for each $v\in V_2(G)$ take a set $M_v$ of $b$ new vertices. Let $H$ be the bigraph with $V_i(H)=\bigcup_{v\in V_i(G)}M_i$
for $i = 1,2$, in which if $u,v$ are $G$-adjacent then $M_u$ is $H$-complete to $M_v$, and if $u,v$ are 
not $G$-adjacent
then $M_u$ is $H$-anticomplete to  $M_v$. We say $H$ is obtained from $G$ by {\em $(a,b)$-multiplication}.
By appropriate multiplication, we can convert an unbalanced parade to a balanced one, and it turns out that all the important properties
of the output of  \ref{usemonotone} are preserved under this. That will allow us to prove:

\begin{thm}\label{usemonotone2}
Let $\delta\ge 2$ and $\eta\ge 0$ be integers.
Let $T$ be a rooted tree bigraph, such that every vertex has degree at most $\delta+1$,
and every path
from root to leaf has length less than $\eta$. Let $\tau=\delta^{\eta+1}$ and
$\lambda = 2^{-30\delta}\delta^{-1-\eta}$.
Let $G$ be a bigraph with a parade
$\mathcal{P}$ of length at least $(K,K)$ where $K= (32\delta+4)\tau+2$, such that
$\mathcal{P}$ is $\lambda$-concave, $(2^{-30\delta},\tau)$-support-invariant and $\tau$-support-uniform.
Let $\mathcal{P}$ have width $(W_1,W_2)$.
If $G$ is $\vare$-coherent where $\vare \le 2^{-30\delta-1}\min(W_1/|V_1(G)|,W_2/|V_2(G)|)$, then there is a $\mathcal{P}$-rainbow copy of $T$.
\end{thm}
\Proof By moving to a sub-parade, we may assume that $\mathcal{P}$ has length $(K,K)$.
Let $|V_i(G)|=n_i$ for $i = 1,2$. Let $H$ be obtained from $G$ by $(W_2,W_1)$-multiplication, with corresponding
sets $M_v\;(v\in V(G))$. Let $\mathcal{P}=(A_1\LL A_K;B_1\LL B_K)$, and let $\mathcal{P}'=(A_1'\LL A_K';B_1'\LL B_K')$, where
each $A_i'=\bigcup_{v\in A_i}M_v$, and $B_i'=\bigcup_{v\in B_i}M_v$ is defined similarly. Thus $\mathcal{P}'$ is a balanced parade in $H$.

$H$ is not yet a balanced bigraph, but we remedy that as follows. Let $U_1=A_1\cup\cdots\cup A_K$, and $U_2=B_1\cup\cdots\cup B_K$.
Thus $|U_2|/|U_1| = W_2/W_1$.
For $i = 1,2$, choose $V_i\subseteq V_i(G)$ with $U_i\subseteq V_i$, such that $|V_2|/|V_1|=|U_2|/|U_1|$, with $V_1$ maximal.
Then $\mathcal{P}$ is a parade in $G[V_1\cup V_2]$, satisfying $|A_i|/|V_1|=|B_j|/|V_2|$ for all $i,j$. Also,
$W_i +|V_i|> n_i$ for some $i\in \{1,2\}$, from the maximality of $V_1$.

Let $G'$ be the sub-bigraph of $H$ induced on the union of the sets $M_v\;(v\in V_1\cup V_2)$. Thus $G'$ is a balanced bigraph, and
$\mathcal{P}'$ is a balanced parade in it.
Let 
$$\vare'= 2\vare \max\left(\frac{W_1/n_1}{W_2/n_2},\frac{W_2/n_2}{W_1/n_1}\right)\ge 2\vare.$$
\\
\\
(1) {\em For $i = 1,2$, $\vare n_i\le \vare' |V_i|$.}
\\
\\
For this we may assume $i=1$, from the symmetry.  Since $\vare'\ge 2\vare $, we may assume that
$n_1>2|V_1|\ge |V_1|+W_1$. 
Since $W_i +|V_i|> n_i$ for some $i\in \{1,2\}$,
it follows that $W_2 +|V_2|> n_2$, and in particular $|V_2|\ge n_2/2$. So
$|V_1|=|V_2|(W_1/W_2)\ge W_1 n_2/(2W_2)$, and so $n_1\le 2\left((W_2/n_2)/(W_1/n_1)\right)|V_1|$.
Hence $\vare n_1 \le \vare'|V_1|$. This proves (1).
\\
\\
(2) {\em $G'$ is $\vare'$-coherent.}
\\
\\
Let $w\in V(G')$, with $w\in M_v$ say. From the symmetry we may assume that $w\in V_1(G)$. Then the set of $G'$-neighbours of $w$
is the union of the sets $M_u$ over all $G$-neighbours $u$ of $v$ with $u\in V_2$. The number
of such $u$ is less than $\vare n_2\le \vare' |V_2|$ by (1); and so $w$ has degree less than $\vare'|V_2(G')|$.

Now suppose that there are subsets $Z_i\subseteq V_i(G')$ of cardinality at least $\vare'|V_i(G')|$ for $i = 1,2$, such that
$Z_1$ is $G'$-anticomplete to $Z_2$. For $i = 1,2$, let $Y_i$ be the set of all $v\in V_i$ such that $M_v\cap Z_i\ne \emptyset$.
It follows that $|Y_1|\ge |Z_1|/W_2$, and $|Y_2|\ge |Z_2|/W_1$, and $Y_1$ is $G$-anticomplete to $Y_2$.
Consequently $|Y_i|< \vare n_i$ for some $i\in \{1,2\}$, and from the symmetry we may assume that $i=1$. Thus
$|Z_1|\le \vare n_1 W_2$, and since $|Z_1|\ge \vare' |V_1(G')|$, it follows that
$ \vare n_1 W_2 > \vare' |V_1(G')|$. But $|V_1(G')|=W_2|V_1|$, and so
$ \vare n_1 W_2 > \vare'W_2|V_1|$, that is, $ \vare n_1 > \vare'|V_1|$, contrary to (1).
This proves (2).
\\
\\
(3) {\em $\mathcal{P}'$ is $\lambda$-concave, $(2^{-30\delta},\tau)$-support-invariant and $\tau$-support-uniform.}
\\
\\
To see that $\mathcal{P}'$ is $\lambda$-concave, let $Y'\subseteq B_1'\cupcup B_L'$. Let $Y$ be the set of all $v\in B_1\cupcup B_L$
such that $Y'\cap M_v\ne \emptyset$. Thus for $1\le i\le K$, the set of vertices in $A_i'$ with a $G'$-neighbour in $Y'$ is
precisely the union of the sets $M_u$ over all $u\in A_i$ that have a $G$-neighbour in $Y$.
So $Y'$ $\lambda$-covers $A_i'$ if and only if $Y$ $\lambda$-covers $A_i$,
and the same for $\lambda$-missing; and so $\mathcal{P}'$ is $\lambda$-concave.

We observe that
for every ordered tree bigraph $S$ with at most $\tau$ vertices,
the trace of $S$ relative to $\mathcal{P}'$ equals its
trace of $S$ relative to $\mathcal{P}$ (because each block of $\mathcal{P}'$ only contains
at most one vertex of the tree).

To see that $\mathcal{P}'$ is $(2^{-30\delta},\tau)$-support-invariant, let $\mathcal{Q}'=(C_1'\LL C_K';D_1'\LL D_K')$
be a contraction of $\mathcal{P}$ such that $|C_i'|\ge 2^{-30\delta}|A_i'|$ and $|D_i'|\ge 2^{-30\delta}|B_i'|$ for $i = 1,2$.
For $1\le i\le K$, there are at least $|C_i'|/W_2$ vertices $v\in A_i$ such that $M_v\cap C_i'\ne \emptyset$.
Let $C_i$ be a set of $\lceil |C_i'|/W_2\rceil$ such vertices. Define $D_1\LL D_k$ similarly. Then
$\mathcal{Q}=(C_1\LL C_K;D_1\LL D_K)$
is a contraction of $\mathcal{P}$, and $|C_i|\ge 2^{-30\delta}|A_i|$ and $|D_i|\ge 2^{-30\delta}|B_i|$ for $1\le i\le K$.
Let $S$ be an ordered tree bigraph with at most $\tau$ vertices. We need to show that for every such $S$, its
trace relative to $\mathcal{Q}'$ equals its trace relative to $\mathcal{P}'$.
The trace of $S$ relative to $\mathcal{Q}'$ is a subset of its trace relative to $\mathcal{P}'$,
and we must show the converse inclusion. Thus, let $(H,K)$ belong to the trace of $S$ relative to $\mathcal{P}'$.
Then $(H,K)$ belongs to the trace of $S$ relative to $\mathcal{P}$, as we saw above.
Since $\mathcal{P}$ is $(2^{-30\delta},\tau)$-support-invariant,
the trace of $S$ relative to $\mathcal{P}$ equals the trace of $S$ relative to $\mathcal{Q}$; and
so $(H,K)$ belongs to the trace of $S$ relative to $\mathcal{Q}$. But then it belongs to the trace relative to $\mathcal{Q}'$.
This proves that $\mathcal{P}'$ is $(2^{-30\delta},\tau)$-support-invariant.

Finally, that $\mathcal{P}'$ is $\tau$-support-uniform is clear, since for every ordered tree bigraph $S$ with at most $\tau$ vertices,
the trace of $S$ relative to $\mathcal{P}'$ equals its
trace relative to $\mathcal{P}$. This proves (3).

\bigskip

Now $\mathcal{P}'$ has width $(W,W)$ where $W=W_1W_2$. From the symmetry we may assume that $W_1/n_1\ge W_2/n_2$.
By hypothesis, $2\vare\le 2^{-30\delta}W_2/n_2$; and since $\vare'=2\vare(W_1/n_1)/(W_2/n_2)$, it follows that
$$\vare'\le 2^{-30\delta}W_1/n_1\le   2^{-30\delta}W_1/|V_1|=2^{-30\delta}W/|V_1(G')|.$$
Since $|V_1(G')|=|V_2(G')|$, 
we deduce, from (2), (3) and \ref{usemonotone} applied to $G', \mathcal{P}'$ and $\vare'$,
that there is a $\mathcal{P}'$-rainbow copy of $T$ in $G'$,
and hence there is a $\mathcal{P}$-rainbow copy of $T$ in $G$. This proves \ref{usemonotone2}.~\bbox

\section{Producing a concave parade}

In this section we apply \ref{usemonotone2} to deduce \ref{mainthm}.

\begin{thm}\label{fixedpoint}
Let $\tau\ge 1$ be an integer, and let $0<\kappa\le 1$.
Let $\mathcal{P}=(A_1\LL A_K; B_1\LL B_L)$ be a parade in a bigraph $G$.
Then there is a contraction $\mathcal{P}'=(A_1'\LL A_K';B_1'\LL B_L')$ of $\mathcal{P}$, such that
\begin{itemize}
\item $|A_i'|\ge \kappa^{2^{K+L}\tau^\tau}|A_i|$ and $|B_j'|\ge  \kappa^{2^{K+L}\tau^\tau}|B_j|$ for $1\le j\le L$; and
\item $\mathcal{P}'$ is $(\kappa,\tau)$-support-invariant.
\end{itemize}
\end{thm}
\Proof
Let $\mathcal{P}=(A_1\LL A_K; B_1\LL B_L)$ be a parade in a bigraph $G$.
We define the {\em trace-cost} of
a contraction $\mathcal{P}'$ of $\mathcal{P}$ to be the sum of the cardinality of the
trace of $T$, summed
over all nonisomorphic ordered tree bigraphs $T$ with at most $\tau$ vertices.
The cardinality of the trace of any given ordered tree bigraph $T$ is at most $2^{K+L}$,
and up to isomorphism there are at most $\tau^\tau$ ordered tree bigraphs $T$ with at most $\tau$ vertices.
Hence the trace-cost
of $\mathcal{P}$ is at most $2^{K+L}\tau^\tau$.

There are integers $t\ge 0$ (for instance $t=0$) such that 
there is a
contraction 
$$\mathcal{P}'=(A_1'\LL A_K';B_1'\LL B_L')$$ of $\mathcal{P}$
with $|A_i'|\ge \kappa^t|A_i|$ for $1\le i\le K$, and $|B_j'|\ge \kappa^t|B_j|$ for $1\le j\le L$, and with trace-cost at most $2^{K+L}\tau^\tau-t$.
Since trace-cost is nonnegative, it follows that every such $t$ satisfies $t\le 2^{K+L}\tau^\tau$, and so we can choose
$t$ maximum with the stated property.
Let $\mathcal{P}''=(A_1''\LL A_K'';B_1''\LL B_L'')$ be a contraction of $\mathcal{P}'$ such that $|A_i''|\ge \kappa|A_i'|$
for $1\le i\le K$,
and $|B_j''|\ge \kappa|B_j'|$
for $1\le j\le L$. For every ordered tree bigraph $T$, the trace of $T$ relative to $\mathcal{P}''$ is a subset of the trace of $T$
relative to $\mathcal{P}$, and so from the choice of $t$, equality holds for every $T$ with at most $\tau$ vertices, that is,
$\mathcal{P}'$ is $(\kappa,\tau)$-support-invariant.
This proves \ref{fixedpoint}.~\bbox

There is a bipartite version of Ramsey's theorem for uniform hypergraphs:
\begin{thm}\label{ramsey}
For all integers $a,b,c$ there exists $N$ with the following property. Let $A,B$ be two disjoint sets both of cardinality at least 
$N$; let $\mathcal{F}$ be the set of all subsets of $A\cup B$ that contain exactly $a$ vertices of $A$ and $b$ vertices of $B$; and let $\mathcal{H}\subseteq \mathcal{F}$.
Then there exists $A'\subseteq A$ and $B'\subseteq B$, with $|A'|=|B'|=c$, such that one of $\mathcal{H}, \mathcal{F}\setminus \mathcal{H}$ contains no subset
of $A'\cup B'$.
\end{thm}

By iterated applications of \ref{ramsey} (one for each ordered tree bigraph $T$
 with at
most $\tau$ vertices)
we deduce:

\begin{thm}\label{getuniform}
Let $k, \tau\ge 0$ be integers; then there exists an integer $K\ge 0$ with the following property. Let
$\mathcal{P}$ be a parade of length $(K,K)$
in a bigraph $G$, and let $0< \lambda\le 1$. Then $\mathcal{P}$
has a sub-parade of length $(k,k)$ which is $\tau$-support-uniform.
\end{thm}
Combining \ref{fixedpoint} and \ref{getuniform}, we obtain:

\begin{thm}\label{getminor}
Let $k,\tau\ge 1$ be integers, and $0<\kappa\le 1$; then there exist an integer $K$ with the following property.
Let $\mathcal{P}$ be a parade of length at least $(K,K)$ and width $(W_1,W_2)$ in a bigraph $G$.
Then there is a minor $\mathcal{P}'$ of $\mathcal{P}$, with length $(k,k)$ and width at least 
$$(\kappa^{2^{2K}\tau^\tau}W_1,\kappa^{2^{2K}\tau^\tau}W_2),$$ 
such that
$\mathcal{P}'$ is $(\kappa,\tau)$-support-invariant and $\tau$-support-uniform.
\end{thm}
\Proof
Let $K$ satisfy \ref{getuniform};
then we claim it satisfies \ref{getminor}. Let
$\mathcal{P}=(A_1\LL A_K; B_1\LL B_K)$ be a parade in a bigraph $G$, of width $(W_1,W_2)$. By \ref{fixedpoint}
there is a contraction $\mathcal{P}'=(A_1\LL A_K;B_1'\LL B_K')$ of $\mathcal{P}$, with width
at least 
$$(\kappa^{2^{2K}\tau^\tau}W_1,\kappa^{2^{2K}\tau^\tau}W_2),$$ 
such that
$\mathcal{P}'$ is $(\kappa,\tau)$-support-invariant.
By \ref{getuniform} applied to $\mathcal{P}'$, the result follows, since being $(\kappa,\tau)$-support-invariant is inherited
by sub-parades. This proves \ref{getminor}.~\bbox

We need the following lemma. 

\begin{thm}\label{getconcave}
Let $K,\ell,t>0$ be integers, and let $0<\kappa\le 1/2$. Let $\lambda=2\kappa$, let $r=\lceil t/\kappa\rceil$, 
and let $L=r\ell$. Let $\mathcal{P}=(A_1\LL A_K;B_1\LL B_L)$ be a parade in a bigraph $G$, of width $(W_1,W_2)$. 
Suppose that 
there do not exist $1<q_0<\cdots< q_{2t}\le L$ and a subset $X\subseteq A_1\cupcup A_K$ such that
$X$ $2\kappa$-covers $B_{q_t}$, and 
$X$ $\kappa$-misses $B_{q_j}$ for all $j\in \{0\LL 2t\}\setminus \{t\}$.
For $1\le j\le \ell$, let
$C_j$ be the union of the sets $B_i$ for all $i$ with $r(j-1)<i\le rj$.
Then the parade $\mathcal{C}=(A_1\LL A_K;C_1\LL C_{\ell})$ is $\lambda$-bottom-concave.
Moreover, for $\tau>0$, if $\mathcal{P}$ is $\tau$-support-uniform then
\begin{itemize} 
\item so is $\mathcal{C}$; 
\item for every ordered tree bigraph $T$ with $|T|\le \tau$ and $|V_2(T)|\le \ell$, if its trace relative to $\mathcal{P}$ is nonempty then its trace relative to $\mathcal{C}$ is nonempty; and
\item if in addition for some $\kappa'\ge \kappa$, $\mathcal{P}$ is 
$(\kappa',\tau)$-support-invariant, then $\mathcal{C}$ is also $(\kappa',\tau)$-support-invariant.
\end{itemize}
\end{thm}
\Proof
Let $1\le h_1<h_2<h_3\le \ell$, and suppose that there exists $X\subseteq A_1\cupcup A_K$, such that
$X$ $\lambda$-covers $C_{h_2}$ and $\lambda$-misses $C_{h_1}$ and $C_{h_3}$. 
Since $X$ $\lambda$-covers $C_{h_2}$, there are at least $\lambda rW_2$ vertices in $C_{h_2}$ with a neighbour in $X$, and so there
exists $r_j$ with $r(h_2-1)< r_j\le rh_2$ such that at least $\lambda W_2$ vertices in $B_{r_j}$ have a neighbour in $X$. Hence
$X$ $\lambda$-covers $B_{r_j}$. From the hypothesis, either $X$ $\kappa$-misses 
$B_i$ for fewer than $t$
values of
$i$ with $r(h_1-1)<i\le rh_1$, or $X$ $\kappa$-misses $B_i$ for fewer than $t$ values of
$i$ with $r(h_3-1)<i\le rh_3$, and from the symmetry we may assume the former. Consequently
there are fewer than
$tW_2+(r-t)\kappa W_2$ vertices in $C_{h_1-1}$ with  no neighbour in $X$. Since $X$ $\lambda$-misses $C_{h_1}$, it follows that
$$tW_2+(r-t)\kappa W_2> \lambda|C_{h_1}|=2\kappa rW_2,$$
and so
$t(1-\kappa)> \kappa r,$
contrary to the choice of $r$. This proves that $\mathcal{C}$ is $\lambda$-bottom-concave. 

If $\mathcal{P}$ is $\tau$-support-uniform then clearly so is $\mathcal{C}$. Let $T$ be an ordered tree bigraph, with $|T|\le \tau$
and $|V_2(T)|\le \ell$ that has nonempty
trace relative to $\mathcal{P}$. Let $|V_1(T)|=s$ and $|V_2(T)|=t$. Thus $t\le \ell$. 
Since $\mathcal{P}$ is $\tau$-support-uniform, the trace of $T$
consists of all pairs $(I,J)$ where $I\subseteq \{1\LL K\}$ with $|I|=s$ and $J\subseteq \{1\LL L\}$ with $|J|=t$. In particular,
$(\{1\LL s\},\{r,2r\LL tr\})$ belongs to the trace of $T$ relative to $\mathcal{P}$; and so $(\{1\LL s\},\{1,2\LL t\})$ belongs to
the trace of $T$ relative to $\mathcal{C}$. This proves the second bullet.

It remains to show that if in addition $\mathcal{P}$ is
$(\kappa',\tau)$-support-invariant then so is $\mathcal{C}$.
Let $T$ be an ordered tree bigraph, with $|T|\le \tau$. Let $|V_1(T)|=s$ and $|V_2(T)|=t$ say.
Let $1\le p_1<\cdots<p_s\le K$ and $1\le q_1<\cdots < q_t\le \ell$, such that there is a 
$(A_{p_1}\LL A_{p_s};C_{q_1}\LL C_{q_t})$-rainbow copy of $T$.
For $1\le i\le s$ let $A_{p_i}'\subseteq A_{p_i}$, and for
$1\le j\le t$ let $C_{q_j}'\subseteq C_{q_j}$, where $(A_{p_1}'\LL A_{p_s}';C_{q_1}'\LL C_{q_t}')$ is a parade with width at least
$(\kappa' W_1, \kappa' rW_2)$. We must show that there is an $(A_{p_1}'\LL A_{p_s}';C_{q_1}'\LL C_{q_t}')$-rainbow copy of $T$.

For $1\le j\le t$, since $|C_{q_j}'|\ge \kappa' rW_2$,
there exists $g_j$ with $r(q_j-1)< g_j\le rq_j$ such that $|C_{q_j}'\cap B'_{g_j}|\ge \kappa' W_2\ge \kappa W_2$.
Choose $D_{q_j}\subseteq C_{q_j}'\cap B'_{g_j}$ of cardinality $\lceil \kappa W_2\rceil$ for $1\le j\le t$.
Since $\mathcal{P}$ is $\tau$-support-uniform and $(\kappa,\tau)$-support-invariant, it follows
that there is a copy of $T$ that is 
$(A_{p_1}'\LL A_{p_s}';D_{q_1}\LL D_{q_s})$-rainbow
and hence $(A_{p_1}'\LL A_{p_s}';C_{q_1}'\LL C_{q_t}')$-rainbow.
This proves~\ref{getconcave}.~\bbox

\begin{thm}\label{rainbow2}
For every tree bigraph $T$, there exist $d>0$ and an integer $K$, such that,
for every bigraph $G$ with a parade $\mathcal{P}$ of length at least $(K,K)$, if for some $\vare>0$, $G$ is $\vare$-coherent 
and $\mathcal{P}$ has width at least $(\vare d |V_1(G)|,\vare d|V_2(G)|)$,
then
there is a $\mathcal{P}$-rainbow copy of $T$ in $G$.
\end{thm}
\Proof
We proceed by induction on $|T|$, and may assume that $|T|\ge 2$. Choose $\delta\ge 2$ and $\eta\ge 0$ such that $T$ is a
sub-bigraph of $T(\delta,\eta)$ and of $\tilde{T}(\delta,\eta)$. (The latter were defined within the proof of \ref{usemonotone}.)
Let $\lambda = 2^{-9\delta-1}\delta^{-1-\eta}$, and $\kappa=\lambda/2$.
Let $r= \lceil (|T|-1)/\kappa\rceil$.
From the inductive hypothesis, 
there exist $K', d'$ such that for  every tree bigraph $T'$ with $|T'|<|T|$, the theorem is satisfied with $T',K',d'$ replacing $T,K,d$.
By increasing $K'$, we may assume that $K'\ge 6r\delta^{\eta+2}$, and $K'>2|T|+1$, and $K'$ is a multiple of $4r$.
Let $\ell=K'/(4r)$.
Let $\tau=\delta^{\eta+1}$.
Let $K$ satisfy \ref{getminor} with $k$ replaced by $K'$.
Let
$$d=\kappa^{-2^K\tau^\tau}\max(d', 2^{9\delta}/r).$$
We claim that $K,d$ satisfy \ref{rainbow2}. Let $\mathcal{P}$
be a parade in a bigraph $G$, of length $(K,K)$ and width $(W_1,W_2)$, where $W_i\ge \vare d |V_i(G)|$ for $i = 1,2$, and 
$G$ is $\vare$-coherent.
We assume (for a contradiction) that there
is no $\mathcal{P}$-rainbow copy of $T$. By \ref{getminor}, there is a minor $\mathcal{P}'$
of $\mathcal{P}$ of length $(K',K')$ and width at least $\kappa^{2^{2K}\tau^\tau}(W_1,W_2)$, such that
$\mathcal{P}'$ is $\tau$-support-uniform and $(\kappa,\tau)$-support-invariant.
Let $\mathcal{P}'=(A_1\LL A_{K'};B_1\LL B_{K'})$, and let its width be $(w_1,w_2)$.
Let $t=|T|$. From the symmetry we may assume that some vertex of $V_1(T)$ has degree one in $T$.
\\
\\
(1) {\em There do not exist $1\le r_0<\cdots< r_{2t}\le K'$, such that for some $X\subseteq A_1\cupcup A_{K'}$, 
$X$ $2\kappa$-covers $B_{r_t}$, and $X$ $\kappa$-misses $B_{r_i}$ for all $i\in \{0\LL 2t\}\setminus \{t\}$.}
\\
\\
Suppose that such $X$ and $r_0\LL r_{2t}$ exist. Let $B_{r_t}'$ be a set of $\lceil \kappa w_2\rceil $ vertices in $B_{r_t}$ that have 
a neighbour in $X$;
for $0\le i\le 2t$ with $i\ne t$, let $B_{r_i}'$ be a set of $\lceil \kappa w_2\rceil $ vertices in $B_{r_i}$ that do not have a neighbour in $X$; and for
$i\in \{1\LL K'\}$ with $i\ne r_0\LL r_{2t}$, let $B_i'$ be a subset of $B_i$ of cardinality $\lceil \kappa w_2\rceil $. 
Partition $\{1\LL K'\}$ into two sets $I_1,I_2$, both of cardinality $K'/2\ge |T|$. Let $X_1$ be the intersection of $X$
with the union of the sets $A_i(i\in I_1)$, and define $X_2$ similarly. At least one of $X_1,X_2$ $\kappa$-covers $B_{r_t}$, since
their union $2\kappa$-covers $B_{r_t}$, and from the symmetry we may assume that $X_2$ $\kappa$-covers $B_{r_t}$.
Since $|I_1|\ge K'/2\ge t$,
Let $J=\{r_0,r_1\LL r_{2t}\}$. Then the parade $\mathcal{P}''=(A_i\;(i\in I_1); B_j'\;(j\in J))$ has width at least
$(w_1, \kappa w_2)$.

Let $v_0\in V_1(T)$ be a vertex of $T$ with degree one, and let $u_0$ be its neighbour.
Let $T'=T\setminus \{v_0\}$.
From the inductive hypothesis, there is a $\mathcal{P}'$-rainbow copy of $T'$, since for $i = 1,2$,
$$\kappa^{2^{2K}\tau^\tau}W_i\ge \kappa^{2^{2K}\tau^\tau}(\vare d |V_i(G)|)\ge d'\vare|V_i(G)|.$$
Hence there is an ordered bigraph $S$, obtained from $T'$ by ordering $V_1(T')$ and $V_2(T')$, with nonempty trace relative to $\mathcal{P}'$. 

Since  $\mathcal{P}'$ is $\tau$-support-uniform, and $(\kappa,\tau)$-support-invariant,
and the trace of $S$  relative to $\mathcal{P}'$
is nonempty, and $|I_1|\ge |V_1(T)|$, there is an isomorphism from $S$ to a $\mathcal{P}''$-rainbow induced sub-bigraph $H$ of $G$, 
with $\mathcal{P}''$-ordering
isomorphic to $S$, where $u_0$ is mapped to a vertex of $H$ in $B'_{r_t}$, say $u$.
Choose $v\in X$ adjacent to $u$; such a vertex exists since PP
$X$ covers $B_{r_t}'$. But then $v$ has no other neighbour in $V(H)$, since $X$ is anticomplete to $B_{r_i}'$
for all $\in I_2\setminus \{t\}$; and $v\in A_h$ where $h\in I_2$.
Thus adding $v$ to $H$ gives a $\mathcal{P}$-rainbow copy of $T$, a contradiction. This proves (1).

\bigskip
We recall that $K'=4r\ell$.
For $1\le i\le 4\ell$, let
$C_i$ be the union of the sets $B_j$ for all $j$ with $r(i-1)<j\le ri$.
Then $\mathcal{C}=(A_1\LL A_{K'}; C_1\LL C_{4\ell})$ is a parade, of width $(w_1,rw_2)$, and $rw_2\ge r\kappa^{2^{2K}\tau^\tau}W_2$.
By \ref{getconcave}, $\mathcal{C}$ is $2\kappa$-bottom-concave, $\tau$-support-uniform and $(\kappa,\tau)$-support-invariant.

Let $1\le i\le K'$, and let $X\subseteq A_i$ be such that $X$ $\lambda$-covers each of $C_1\LL C_{4\ell}$. This is possible, because
$A_i$ $\lambda$-covers each of $C_1\LL C_{4\ell}$, since $\lambda\le 1/2$
and $w_1\ge \vare|V_1(G)|$, and $rw_2\ge \vare|V_2(G)|$. Let us say $j\in \{1\LL 4\ell\}$ is a {\em pit} for $X$ if $X$ $\lambda$-misses $C_j$.
Since $\mathcal{C}$ is $\lambda$-bottom-concave, there are at most two pits for $X$,
and if there are two then they are consecutive. On the other hand, if we choose $X$ minimal, then 
it follows by deleting one vertex from $X$ that for some $j\in \{1\LL 4\ell\}$, fewer than $\lambda |C_j|+\vare|V_2(G)|$ vertices in $C_j$
have a neighbour in $X$, and since $\lambda |C_j|+\vare|V_2(G)|\le (1-\lambda)|C_j|$, it follows that $j$ is a pit for $X$.
Thus there is either exactly one pit, or exactly two consecutive pits for $X$. Let us say $A_i$ is {\em left-first} if
there is a choice of $X\subseteq A_i$ such that $X$ $\lambda$-covers each of $C_1\LL C_{4\ell}$, and $X$ has a pit $j$ with $j>2\ell$;
and {\em right-first} if
there is a choice of $X\subseteq A_i$ such that $X$ $\lambda$-covers each of $C_1\LL C_{4\ell}$, and $X$ has a pit $j$ with $j\le 2\ell$.
Each $A_i$ is therefore either left-first or right-first, possibly both; and so we may assume (by reversing the order of 
$C_1\LL C_\ell$ if necessary) that at least half of them are left-first.
Let $I\subseteq \{1\LL K'\}$ with $|I|=K'/2=2r\ell$ such that $A_i$ is left-first for each $i\in I$.

Let $i\in I$, and choose $X\subseteq A_i$ such that $X$ $\lambda$-covers each of $C_1\LL C_{4\ell}$, and $X$ has a pit $j$ with $j>2\ell$.
Choose $X^{2\ell}\subseteq X$ minimal such that $X^{2\ell}$ $\lambda$-covers $C_{2\ell}$, and for $2\ell> i\ge 1$ in turn, 
inductively choose $X^i\subseteq X^{i+1}$ minimal such that $X^i$ $\lambda$-covers $C_i$. To show that this is possible, we will
prove inductively that for $2\ell\ge  i\ge 1$:
\begin{itemize}
\item $X^i$ does not $\lambda$-miss any of $C_1\LL C_{i-1}$; 
\item $X^i$ both $\lambda$-covers and $\lambda$-misses $C_i$; 
\item $X^i$ does not $\lambda$-cover any of $C_{i+1}\LL C_{2\ell}$; and
\item fewer than $\lambda |C_i|+\vare|V_2(G)|$ vertices in $C_i$
have a neighbour in $X^i$.
\end{itemize}
Suppose then that either $i=2\ell$ or $X^{i+1}$ satisfies the four bullets; then the choice of $X^i$ is possible, and it remains to show that $X^i$
satisfies the four bullets. Certainly $X^i$ $\lambda$-covers $C_i$ from its definition; and from the minimality of $X^i$, 
fewer than $\lambda |C_i|+\vare|V_2(G)|\le (1-\lambda)|C_i|$ vertices in $C_i$
have a neighbour in $X^i$, and so $X^i$ $\lambda$-misses $C_i$. Since it also $\lambda$-misses $C_j$, because $j$ is a pit for $X$ and
$X^i\subseteq X$, and $j>2\ell$, we deduce from concavity that $X^i$ does not $\lambda$-cover any of $C_{i+1}\LL C_{2\ell}$. Thus
the second, third and fourth bullets hold. But also, since $X^i$ $\lambda$-covers $C_i$ and $\lambda$-misses $C_j$, concavity 
implies that $X^i$ does not $\lambda$-miss any of $C_1\LL C_{i-1}$. Thus all four bullets hold. This completes the 
inductive definition of $X^1\LL X^{2\ell}$. Let us call the sequence $(X^1\LL X^{2\ell})$ a {\em $\lambda$-ladder} in $A_i$.

For each $i\in I$, choose a $\lambda$-ladder $(X_i^1\LL X_i^{2\ell})$ in $A_i$. For $1\le j\le 2\ell$, let us say that
$v\in C_j$ is {\em unwanted} if for some $i\in I$,  either $v$ has a neighbour in $X_i^j$, or $v$ has no neighbour in $X_i^{j'}$ 
for some $j'$ with $j<j'\le 2\ell$ (and hence $j<2\ell$, and $v$ has no neighbour in $X_i^{j+1}$). We say $v$ is {\em wanted} if it 
is not unwanted. (Note that if $v\in C_j$ is wanted, then for each $i\in I$, $v$ has no neighbour in any of $X^i_1\LL X^i_j$.)
It follows that the total number of unwanted vertices in $C_j$ is at most $|I|(2\lambda |C_j|+\vare|V_2(G)|)$,
since for each $i\in I$, there are at most $\lambda |C_j|+\vare|V_2(G)|$ vertices in $C_j$ with a neighbour in $X^j_i$, and
at most $\lambda |C_j|$ vertices in $C_j$ that have no neighbour in  $X_i^{j+1}$ (when $j<2\ell$). Since
$|I|(2\lambda |C_j|+\vare|V_2(G)|)\le |C_j|/2$
it follows that at least $|C_j|/2$ vertices in $C_j$ are wanted. 
For $1\le j\le 2\ell$, let $C_j'\subseteq C_j$ be a set of vertices that are wanted, with $|C_j'|=\lceil |C_j|/2\rceil$.
\\
\\
(2) {\em There do not exist $q_0\LL q_{2t}\in I$ with $q_0<q_1<\cdots <q_{2t}$, such that for some $Y\subseteq C_1'\cupcup C_{\ell}'$,
$Y$ $2\kappa$-covers $A_{q_t}$, and $Y$ $\kappa$-misses $A_{q_i}$ for all $i\in \{0\LL 2t\}\setminus \{t\}$.}
\\
\\
Suppose that such $q_0\LL q_{2t}, j, Y$ exist. Choose $u\in V_2(T)$ such that all its $T$-neighbours except possibly one have 
degree one in $T$. (This is possible, because if $V_2(T)$ contains a leaf of $T$, let $u$ be that leaf, and if all leaves
belong to $V_1(T)$, let $u$ be a leaf of the tree obtained from $T$ by deleting all leaves.) Let $v$ be a neighbour of $u$
such that all its other neighbours are leaves; and let $u$ have $s$ neighbours different from $v$.
Let $T'$ be obtained from $T$
by deleting $u$ and all its neighbours except $v$. 
Let $A_{q_t}'$ be a subset of $\lceil \kappa W_1\rceil$ vertices in $A_{q_t}$ that have a neighbour in $Y$, and for $0\le i\le 2t$
with $i\ne t$, let $A_{q_i}'$ be a subset of $\lceil \kappa W_1\rceil$ vertices in $A_{q_i}$ that have no neighbour in $Y$.
Let $I_1=\{q_0,q_1\LL q_{2t}\}$.

As in the proof of (1), from the inductive hypothesis, there is a $\mathcal{P}'$-rainbow copy of $T'$.
Hence there is an ordered bigraph $S$, obtained from $T'$ by ordering $V_1(T')$ and $V_2(T')$, with nonempty trace 
relative to $\mathcal{P}'$, and so with nonempty trace relative to $\mathcal{C}$, by the second bullet of \ref{getconcave}. 

Let $\mathcal{C}'=(A_i'\;(i\in I_1);C_j'\;(\ell+1\le j\le 2\ell))$. 
Since  $\mathcal{C}$ is $\tau$-support-uniform, and $(\kappa,\tau)$-support-invariant, and $\kappa\le 1/2$ (and so each $|C_j'|\ge \kappa|C_j|$), 
and the trace of $S$  relative to $\mathcal{C}$
is nonempty, it follows that there is a $\mathcal{C}'$-rainbow induced sub-bigraph $H$ of $G$, with $\mathcal{P}$-ordering
isomorphic to $S$,
where some vertex $v'\in V(H)\cap A_{q_t}$ is mapped by the isomorphism to $v$.
Choose $u'\in Y$ adjacent to $v'$; such a vertex exists since
$Y$ covers $A_{q_t}'$. But then $u'$ has no other neighbour in $V(H)$, since $Y$ misses $A_{q_i}'$
for all $i\in \{0\LL 2t\}\setminus \{t\}$.
Choose $p_1\LL p_s\in I\setminus \{q_0\LL q_{2t}\}$. (This is possible since $|I|\ge 2t+1+s$.) 
Since $u'\in Y\subseteq C_1'\cupcup C_{\ell}'$, there exists $j$ with $1\le j\le \ell$ such that $u'\in C_j'$.
For $1\le i\le s$, choose $x_i\in X^{p_i}_j$
adjacent to $u'$. (This is possible since $u'$ is wanted.)
Then $x_i$ is nonadjacent to all vertices in $V_2(H)$, since all of these vertices are wanted and belong to 
$C_{\ell+1}'\cup\cdots\cup C_{2\ell}'$.
Thus adding $u'$ and $x_1\LL x_s$ to $H$ gives a $\mathcal{P}$-rainbow copy of $T$, a contradiction. This proves (2).

\bigskip

Since $|I|=2r\ell$, we may partition $I$ into $\ell$ ``intervals'' each containing $2r$ elements of $I$, say $I_1\LL I_{\ell}$.
More exactly, we partition $I$ into $I_1\LL I_{\ell}$, where each of these sets has cardinality $2r$, and for $1\le h<i\le \ell$, every element
of $I_h$ is less than every element of $I_i$. For $1\le h\le \ell$ let $D_h$ be the union of the sets $A_i\;(i\in I_h)$.
Then $\mathcal{D}=(D_1\LL D_{\ell}; C_{1}'\LL C_{\ell}')$ is a parade. 
By (2) and \ref{getconcave}, $\mathcal{D}$ is $\lambda$-top-concave. Now
$\mathcal{C}$ is
$\tau$-support-uniform, and $(\kappa,\tau)$-support-invariant, and $\lambda$-bottom-concave; and so
$\mathcal{C}'=(A_i\;(i\in I); C_{1}'\LL C_{\ell}')$ is
$\tau$-support-uniform, $(2\kappa,\tau)$-support-invariant, and $2\lambda$-bottom-concave, since 
$|C_i'|\ge |C_i|/2$ for each $i$. Thus $\mathcal{D}$ is $2\lambda$-bottom-concave; and by \ref{getconcave}, $\mathcal{D}$ is also $\tau$-support-uniform and $(2\kappa,\tau)$-support-invariant.
In summary then, $\mathcal{D}$ is $\tau$-support-uniform, $(2\kappa,\tau)$-support-invariant and $2\lambda$-concave.
Hence $\mathcal{D}$ satisfies the hypotheses of \ref{usemonotone2}, and so there is a $\mathcal{D}$-rainbow, and hence $\mathcal{P}$-rainbow,
copy of $T$, a contradiction. This proves \ref{rainbow2}.~\bbox

Finally we can prove \ref{mainthm}, which we restate:

\begin{thm}\label{mainthm2}
For every forest bigraph $H$, there exists $\vare>0$, such that every $\vare$-coherent bigraph contains $H$.
\end{thm}
\Proof
We may assume that $H$ is a tree bigraph. Let $K,d$ satisfy \ref{rainbow2}. 
We may assume by increasing $d$ that $d\ge 1$. Choose $\vare>0$ such that $2K\vare d\le 1$. 
We claim that every $\vare$-coherent bigraph contains $T$.
Let $G$ be an $\vare$-coherent bigraph. It follows that $|V_1(G)|, |V_2(G)|\ge \vare^{-1}\ge 2Kd\ge 2K$. Hence
for $i = 1,2$, $|V_i(G)|/K\ge \lceil |V_i(G)|/(2K)\rceil$, and so we may choose $K$ subsets of $V_i(G)$, pairwise disjoint and
each of cardinality $\lceil |V_i(G)|/(2K)\rceil \ge \vare d |V_i(G)|$. These sets, in order, form a parade of length $(K,K)$
and width at least $\vare d (|V_1(G)|,|V_2(G)|)$,
and so by \ref{rainbow2}, $G$ contains $H$.
This proves \ref{mainthm2}.~\bbox

\end{document}